\documentclass[12pt, reqno]{amsart}
\usepackage{amsthm}
\usepackage{amssymb}
\usepackage{amsmath}

\def\1{\hbox{1\kern-.35em\hbox{1}}}

\newtheorem{theorem}{Theorem}[section]
\newtheorem*{theorem*}{Theorem}
\newtheorem{lemma}{Lemma}[section]
\newtheorem{proposition}{Proposition}[section]
\newtheorem*{proposition*}{Proposition}
\newtheorem{corollary}{Corollary}[section]
\newtheorem{definition}{Definition}[section]
\newtheorem{remark}{Remark}[section]

\numberwithin{equation}{section}

\newcommand{\bea}{\begin{eqnarray}}
\newcommand{\eea}{\end{eqnarray}}
\newcommand{\be}{\begin{eqnarray*}}
\newcommand{\ee}{\end{eqnarray*}}

\newcommand{\Z}{{\mathbb Z}}

\newcommand{\R}{{\mathbb R}}
\newcommand{\C}{{\mathbb C}}
\newcommand{\Cq}{{\mathbb C(q)}}
\newcommand{\Rq}{{{\mathbb R}(q)}}

\newcommand{\fb}{{\mathfrak b}}
\newcommand{\fg}{{\mathfrak g}}
\newcommand{\fh}{{\mathfrak h}}
\newcommand{\fj}{{\mathfrak j}}
\newcommand{\fk}{{\mathfrak k}}
\newcommand{\fl}{{\mathfrak l}}

\newcommand{\fp}{{\mathfrak p}}
\newcommand{\fq}{{\mathfrak q}}
\newcommand{\fr}{{\mathfrak r}}

\newcommand{\fu}{{\mathfrak u}}

\newcommand{\ad}{{\rm ad}}
\newcommand{\id}{{\rm id}}
\newcommand{\Hom}{{\rm Hom}}
\newcommand{\I}{{\rm I}}
\newcommand{\Uc}{{\Ug^\circ}}
\newcommand{\U}{{\rm U}}

\newcommand{\Ug}{{{\rm U}_q(\fg)}}
\newcommand{\Uh}{{{\rm U}_q(\fh)}}
\newcommand{\Ub}{{{\rm U}_q(\fb)}}
\newcommand{\Ubb}{{{\rm U}_q(\bar\fb)}}
\newcommand{\Ur}{{{\rm U}_q(\fr)}}
\newcommand{\Up}{{{\rm U}_q(\fp)}}
\newcommand{\Ul}{{{\rm U}_q(\fl)}}
\newcommand{\Uq}{{{\rm U}_q(\fq)}}
\newcommand{\Uk}{{{\rm U}_q(\fk)}}
\newcommand{\Uj}{{{\rm U}_q(\fj)}}

\newcommand{\URg}{{{\rm U}_q^{\R}(\fg)}}

\newcommand{\URl}{{{\rm U}_q^{\R}(\fl)}}

\newcommand{\Ag}{{{\mathcal A}(\fg)}}

\newcommand{\cA}{{\mathcal A}}
\newcommand{\cC}{{\mathcal C}}
\newcommand{\cF}{{\mathcal F}}
\newcommand{\cS}{{\mathcal S}}

\newcommand{\tL}{{\tilde L}}
\newcommand{\tR}{{\tilde R}}

\begin{document}

\title[Quantum supergroups and non-commutative bundles]
{Quantum superalgebra representations on cohomology groups of
non-commutative bundles}
\author[R. B. ZHANG]{R. B. ZHANG}
\address{School of Mathematics and Statistics,
University of Sydney, NSW 2006, Australia}
\email{rzhang@maths.usyd.edu.au}

\begin{abstract} Quantum homogeneous supervector bundles
arising from the quantum general linear supergoup are studied. The space of
holomorphic sections is promoted to a left exact covariant functor
from a category of modules over a quantum parabolic sub-supergroup
to the category of locally finite modules of the quantum general
linear supergroup. The right derived functors of this functor
provides a form of Dolbeault cohomology for quantum homogeneous
supervector bundles. We explicitly compute the
cohomology groups, which are given in terms of well understood
modules over the quantized universal enveloping algebra of the
general linear superalgebra.
\end{abstract}

\maketitle
%\tableofcontents

\section{Introduction}
We follow the general philosophy of non-commutative
geometry \cite{Co, Ma88} to study
quantum homogeneous supervector bundles arising
from the quantum general linear supergroup. Our starting
point is the quantized universal enveloping algebra
$\U_q({\mathfrak{gl}}_{m|n})$ (see e.g., \cite{Z93, Zou})
of the complex general linear
superalgebra ${\mathfrak{gl}_{m|n}}$ \cite{Kac, Sc}. As is well
known, $\U_q({\mathfrak{gl}}_{m|n})$ has the structure of a Hopf
superalgebra \cite{Mo, MM}. Thus its dual superspace
$\U_q({\mathfrak{gl}}_{m|n})^*$ acquires a natural associative
superalgebraic structure. The subspace $\cA(\mathfrak{gl}_{m|n})$
of $\U_q({\mathfrak{gl}}_{m|n})^*$ spanned by all the matrix
elements of the finite dimensional $\U_q({\mathfrak{gl}}_{m|n})$-representations
with integral weights forms a Hopf superalgebra, which  may be
considered as the superalgebra of
functions on the quantum general linear supergroup. This Hopf
superalgebra is closely related to the multi-parameter quantization
of the general linear supergroup of \cite{Ma89}, and obviously
contains the Hopf superalgebra $G_q$ of \cite{Z98} as a Hopf
sub-superalgebra.

For any given reductive quantum sub-superalgebra $\Ul$ of
$\U_q({\mathfrak{gl}}_{m|n})$, we consider the subspace
$\cA(\mathfrak{gl}_{m|n}, \fl)$ of $\cA(\mathfrak{gl}_{m|n})$
invariant with respect to left translations under $\Ul$. This
subspace forms a sub-superalgebra of $\cA(\mathfrak{gl}_{m|n})$.
In the spirit of non-commutative geometry
\cite{Co, Ma88}, we shall regard $\cA(\mathfrak{gl}_{m|n}, \fl)$ as (the
superalgebra of functions on) a quantum homogeneous superspace,
and finite type projective $\cA(\mathfrak{gl}_{m|n}, \fl)$-modules
as (spaces of global sections of) quantum supervector bundles on
the quantum homogeneous superspace. As we shall see, such
$\cA(\mathfrak{gl}_{m|n}, \fl)$-modules also admit a natural
$\U_q({\mathfrak{gl}}_{m|n})$-action. This provides an
interesting link between the non-commutative geometry of the
quantum supervector bundles and the representation theory of
$\U_q({\mathfrak{gl}}_{m|n})$. In the context of classical Lie
groups, such a link is very well known and constitutes the subject
of the Bott-Borel-Weil theory \cite{Bo} (see \cite{KV} for a review on
the algebraic theory). For Lie supergroups in
the classical setting, the subject was studied by Penkov \cite{Pe}
(see also \cite{San}).

We shall study a version of Dolbeault cohomology for quantum
homogeneous supervector bundles. Our method is based on
Zuckerman's algebraic theory of induced representations
\cite{EW, KV}, which was generalized to quantum groups by Andersen
and co-workers in \cite{APW}, and also to Lie superalgebras in \cite{San}.
We shall make essential use of results of \cite{APW}.

We promote the space $\cS$ of global sections and
the space $\Gamma$ of holomorphic sections of a quantum
homogeneous supervector bundle to covariant
functors from appropriate categories of modules over quantum
sub-superalgebras of $\U_q({\mathfrak{gl}}_{m|n})$ to the category
of locally finite $\U_q({\mathfrak{gl}}_{m|n})$-modules. The
resultant functors are closely related to a quantum analogue of
the Zuckerman functor, which will be introduced in Section
\ref{Functors}. The derived functors of the `holomorphic section
functor' arising from $\Gamma$ are the Dolbeault cohomology groups
for quantum homogeneous supervector bundles which we seek for.
When the quantum homogeneous supervector bundle $\cS$ is induced
by a finite dimensional irreducible module over a purely even
reductive quantum subalgebra, or a finite dimensional dual Kac
module (defined by \eqref{dKac}) over a general reductive quantum
sub-superalgebra,   we explicitly compute the
cohomology groups in Theorem \ref{main} and Theorem
\ref{general}. The results resemble the classical
Bott-Borel-Weil theorem \cite{Bo} in that the cohomology is
concentrated at one degree. However, the non-trivial cohomology
groups are isomorphic to dual Kac modules over
$\U_q({\mathfrak{gl}}_{m|n})$, thus are not irreducible. The weight
spaces of the cohomology groups are worked out.

The arrangement of the paper is as follows. Section
\ref{Preliminaries} contains some background material on the
general linear superalgebra and its quantized universal enveloping
algebra. Section \ref{Bundles} introduces the notions of quantum
homogeneous superspaces and quantum homogeneous supervector
bundles on them, and study some basic properties of theirs.
Section \ref{Functors} studies induction functors. The Zuckerman
functor is defined, and a version of Frobenius reciprocity is
proven. Properties of various functors and module categories are
also established, which are to be used in the next section. In
Section \ref{BBW} we first present the formulation Dolbeault
cohomology for quantum homogeneous supervector bundles using
results of Section \ref{Functors}, then compute the Dolbeault
cohomology groups for bundles of interest.

\section{Preliminaries} \label{Preliminaries}
This section presents some preliminary material on the general
linear superalgebra and its quantized universal enveloping
algebra, which sets up the stage for the remainder of the paper.
The section also serves to fix notations and conventions.
\subsection{The general linear superalgebra ${\mathfrak{gl}}_{m|n}$}
Through out the paper we shall denote by $\fg$ the complex general
linear superalgebra ${\mathfrak{gl}}_{m|n}$ \cite{Kac, Sc}. Let ${\bf I}=\{1, 2,
\cdots, m+n\}$, and ${\bf I}'=\{1, 2, \cdots, m+n-1\}$. We shall
identify $\fg$ with the Lie superalgebra of $(m+n)\times
(m+n)$-matrices, where the $\Z_2$-grading is specified as follows.
Denote by $e_{a b}$, $a, b\in{\bf I}$, the
$(m+n)\times(m+n)$-matrix unit with all entries being zero except
that at the $(a, b)$ position which is $1$. We declare $e_{a b}$
to be odd if $a\le m <b$ or $a>m\ge b$, and even otherwise.  Then
$\{e_{a b} | a, b\in{\bf I}\}$ forms a homogeneous basis of $\fg$.
The maximal even subalgebra of $\fg$ will be denoted by $\fg_0$,
which is equal to ${\mathfrak{gl}}_m\oplus {\mathfrak{gl}}_n$. Let
$\fg_{+1}=\sum_{i\le m<r}\C e_{i r}$, and $\fg_{-1}=\sum_{i\le
m<r}\C e_{r i}$. Then the odd subspace of $\fg$ is
$\fg_{+1}\oplus\fg_{-1}$.

We fix the Borel subalgebra $\fb$ consisting of the upper
triangular matrices, and take  ${\fh}=\bigoplus_a \C E_{a a}$ as
the Cartan subalgebra.  Let $\{\epsilon_a \,|\, a\in{\bf I}\}$ be
the basis of ${\mathfrak h}^*$ such that $\epsilon_a(E_{b
b})=\delta_{a b}$. The space ${\mathfrak h}^*$ is equipped with a
bilinear form $(\  , \ ): {\mathfrak h}^*\times {\mathfrak h}^*
\rightarrow \C$ such that $(\epsilon_a, \epsilon_b)=\left\{
\begin{array}{l l}\delta_{a b}, &a\le m, \\
-\delta_{a b}, &a> m. \end{array}\right. $ We shall denote by
$\fh^*_{\Z}$ the $\Z$-span of the $\epsilon_a$. The set of roots
of $\fg$ is $\{ \epsilon_a-\epsilon_b | a\ne b\}$, with
$\epsilon_a-\epsilon_b$ being called odd if $a\le m<b$ or $b\le
m<a$, and even otherwise. The set of the positive roots relative
to the Borel subalgebra $\fb$ is $\{ \epsilon_a-\epsilon_b \,|\,
a< b\}$, and the set of simple roots is
$\{\epsilon_a-\epsilon_{a+1} \,|\, a\in{\bf I}' \}$. An element
$\lambda\in\fh^*$ is called dominant if ${2(\lambda, \,
\alpha)}/{(\alpha,\, \alpha)} \in\Z_+ $, for all positive even
roots of $\fg.$ Denote by $2\rho$ the signed-sum of the positive
roots of $\fg$. A $\lambda\in\fh^*$ is called $\fg$-regular if
$(\lambda+\rho, \, \alpha)\ne 0$ for all even roots of $\fg.$

The elements of the following set $\{e_{a, a+1},  e_{a+1, a} |
a\in{\bf I}'\}\cup\{e_{b b} | b\in{\bf I}\}$ generate $\fg$. We
shall call a Lie sub-superalgebra $\fr$ of $\fg$ regular if there
exist subsets $\Theta_{\pm }$ of ${\bf I}'$ and a subset
$\Theta_0$ of ${\bf I}$ such that $\fr$ is generated by the
elements of the set $\{e_{a a} | a\in\Theta_0\}\cup\{e_{b, b+1} |
b\in\Theta_+\}\cup\{e_{c+1, c} |c\in\Theta_-\}$. The Lie
sub-superalgebra $\fr$ is called reductive if $\Theta_0={\bf I}$
and $\Theta_+=\Theta_-$, and is called parabolic if $\Theta_0={\bf
I}$, and either $\Theta_+$ or $\Theta_-$ is equal to ${\bf I}'$.
If $\fr$ is a parabolic Lie sub-superalgebra, then it contains the
reductive Lie sub-superalgebra, called the Levi factor of $\fr$,
generated by the elements of the set $\{e_{a a} |
a\in\Theta_0\}\cup\{e_{b, b+1}, e_{b+1, b} |
b\in\Theta_+\cap\Theta_-\}$. Note that a parabolic Lie
sub-superalgebra of $\fg$ necessarily contains $\fb$ or  the
opposite Borel subalgebra $\bar\fb$ spanned by the lower
triangular matrices.

\subsection{The quantized universal enveloping algebra of ${\mathfrak{gl}}_{m|n}$}
  Let $q$ be an indeterminate, and
denote by $\Cq$ the field of complex rational functions in $q$.
Set $q_a=\left\{\begin{array}{l l} q, & a\le m,\\
q^{-1}, & a>m. \end{array} \right.$ We define the quantized
universal enveloping algebra $\Ug$ \cite{Z93, Zou, Z98} of the general linear
superalgebra $\fg$ to be the unital associative superalgebra over
$\Cq$ with the set of generators \be \{E_{a, a+1}, \, E_{a+1, a}
\, |\, a\in{\bf I}'\}\cup \{K_b, \, K_b^{-1}\, | \, b \in{\bf I}\}
\ee subject to the following relations
\bea
K_a K_a^{-1}=1,
& & K_a^{\pm  1} K_b^{\pm 1} = K_b^{\pm 1}  K_a^{\pm 1}, \nonumber \\
K_a E_{b\ b\pm 1} K_a^{-1} &=&
q^{(\epsilon_a, \, \epsilon_b -\epsilon_{b\pm 1})} E_{b\  b\pm 1}, \nonumber \\
{}[E_{a\,  a+1},\,E_{b+1\,  b}\}& =& \delta_{a b}
\frac{K_a K_{a+1}^{-1} - K_a^{-1} K_{a+1}}{q_a - q_a^{-1}},\nonumber \\
(E_{m\, m+1})^2 &=& (E_{m+1\, m})^2 = 0, \nonumber \\
E_{a\,  a+1} E_{b\,  b+1} &=& E_{b\,  b+1} E_{a\,  a+1},\nonumber \\
E_{a+1\, a} E_{b+1\, b} &=&E_{b+1\, b} E_{a+1\, a}, \ \ \
\vert a - b\vert \ge 2, \nonumber \\
{\mathcal S}^{(+)}_{a \ a\pm 1}&=&{\mathcal S}^{(-)}_{a \ a\pm 1}=0,
\ \ \ a\ne m,\nonumber \\
\{ E_{m-1\, m+2},\ E_{m\, m+1}\} &=&
\{ E_{m+2\, m-1},\ E_{m+1\, m}\} = 0,   \label{quantum}\eea
where $[E_{a\,  a+1},\,E_{b+1\,  b}\}:=E_{a\,  a+1} \, E_{b+1\,  b}
- (-1)^{\delta_{a m}\delta_{b m}} E_{b+1\,  b}\, E_{a\,  a+1}.$
The $E_{m-1\, m+2}$ and $E_{m+2\, m-1}$ are the $a=m-1$, $b=m+1$,
cases of the elements defined by \eqref{generators}, and
\be
{\mathcal S}^{(+)}_{a \ a\pm 1}&=&
(E_{a\, a+1})^2  E_{a\pm 1\, a+1\pm 1} - (q +
q^{-1}) E_{a\, a+1} \ E_{a\pm 1\, a+1\pm 1} \ E_{a\, a+1}\\
& +& E_{a\pm 1\, a+1\pm1 }\ (E_{a\, a+1})^2,    \\
{\mathcal S}^{(-)}_{a \ a\pm 1}&=&
(E_{a+1\, a})^2\,E_{a+1\pm 1\, a\pm 1} - (q +
q^{-1}) E_{a+1\, a}\ E_{a+1\pm 1\, a\pm 1} \ E_{a+1\, a}\\
&+& E_{a+1\pm 1\, a\pm 1}\ (E_{a+1\, a})^2.
\ee
The $\Z_2$ grading of the superalgebra is defined by declaring
the elements $K_a^{\pm 1}$, $\forall a\in {\bf I}$,
and $E_{b\,  b+1}$, $E_{b+1\,  b}$,  $b\ne m$, to be even and
$E_{m\, m+1}$ and $E_{m+1\, m}$ to be odd. Throughout the paper,
we use $[f]$ to denote the parity of the element $f$ of any
$\Z_2$-graded space.

It is well known that $\Ug$
has the structure of a Hopf superalgebra \cite{MM, Mo},
with a co - multiplication \be
\Delta(E_{a\, a+1}) &=& E_{a\,  a+1} \otimes
K_a K_{a+1}^{-1} + 1 \otimes E_{a\, a+1}, \\
\Delta(E_{a+1\, a}) &=& E_{a+1\, a }\otimes 1 + K_a^{-1} K_{a+1}
\otimes E_{a+1\, a}, \\
\Delta(K_a^{\pm 1}) &=&K_a^{\pm 1}\otimes K_a^{\pm 1},
\ee
co - unit
\be
\epsilon(E_{a\, a+1})&=&E_{a+1\, a}=0, \ \ \forall a\in{\bf I}', \\
\epsilon(K_b^{\pm 1})&=&1,  \ \ \ \forall b\in{\bf I},
\ee
and antipode
\be
S(E_{a\, a+1}) &=& - E_{a\, a+1} K_a^{-1} K_{a+1}, \\
S(E_{a+1\, a}) &=& - K_a K_{a+1}^{-1}E_{a+1\, a}, \\
S(K_a^{\pm 1}) &=&K_a^{\mp 1}\otimes K_a^{\mp 1}.
\ee

Let $E_{a b}$, $E_{b a}$, $a<b$, be elements of $\Ug$ defined
by  \bea E_{a b} &=& E_{a c} E_{c b} - q_c^{-1} E_{c b} E_{a c},
\quad a<c<b, \nonumber \\
E_{b a} &=& E_{b c} E_{c a} - q_c E_{c a} E_{b c}, \quad a<c<b.
\label{generators}
\eea The definition is independent of the $c$ chosen \cite{Z93}.
These elements  are the generalization to
$\U_q({\mathfrak{gl}}_{m|n})$ of a similar set of elements for
$U_q({\mathfrak{gl}}_n)$ constructed by Jimbo \cite{Ji}.
They were used in \cite{KT} for the construction of the universal
$R$-matrix of $\U_q({\mathfrak{gl}}_{m|n})$. We mention that the $E_{a b}$
behave very much like the images of $e_{a b}$ in the universal
enveloping algebra of $\mathfrak{gl}_{m|n}$. For example,  $E_{a
b}^2=E_{b a}^2=0$, if $a\le m<b$. More importantly, we have the
following Poincar\'{e}-Birkhoff-Witt theorem for $\Ug$.
\begin{theorem}\cite{Z93, Zou} \label{PBW}
The ordered products of non-negative powers of all the $E_{a b}$,
$a\ne b$, and integer powers of all $K_c$, $c\in{\bf I}$,  with
respect to any linear ordering of the elements of $\{E_{a b}, E_{b
a} |a<b\}\cup\{K_c |c\in{\bf I}\}$ form a basis of $\Ug$.
\end{theorem}

We shall assume that every
$\Ug$-module to be considered in this paper is $\Z_2$-graded.
Let $V$ be a $\Ug$-module. A weight vector $v\in V$ is the
simultaneous eigenvector of all the $K_a$, $a\in{\bf I}$. We shall
say that $v$ is an integral weight vector with weight
$\mu\in\fh^*_\Z$ if \be K_a v = q^{(\mu, \epsilon_a)} v, & &
\forall a.\ee We shall denote by $L_\lambda$,
$\lambda\in\fh_\Z^*$, the irreducible $\Ug$-module with a highest
weight vector which is an integral weight vector with weight
$\lambda$. The $\lambda$ will be referred to as the highest weight
of $L_\lambda$. It was shown in \cite{Z93} that $L_\lambda$ is
finite dimensional if and only if its highest weight is dominant.

We shall call a Hopf sub-superalgebra of $\Ug$ a quantum
sub-superalgebra. Corresponding to a regular Lie sub-superalgebra
$\fr$ of $\fg$ specified by the sets $\Theta_0$ and
$\Theta_{\pm}$, there exists an associated quantum
sub-superalgebra $\Ur$ generated by the elements of the following
set $\{K_a, \, K_a^{-1} | a\in\Theta_0\}\cup\{E_{b, b+1} |
b\in\Theta_+\}\cup\{E_{c+1, c} |c\in\Theta_-\}$. Important quantum
sub-superalgebras are $\Uh$ and the two quantum Borel
sub-superalgebras $\Ub$ and $\Ubb$. If $\fr$ is a parabolic
(respectively reductive) Lie sub-superalgebra of $\fg$, then $\Ur$
will be called a parabolic (respectively reductive) quantum
sub-superalgebra of $\Ug$. If $\Up$ is parabolic with the Levi
factor $\Ul$, then we have the Hopf superalgebra inclusions
$\Ul\subset\Up\subset\Ug$.

For later use, we consider here a particular module over a
reductive quantum sub-superalgebra $\Ul$. Let $\rm{U}_q(\fl_{0})$
$=$ $\Ul\cap\rm{U}_q(\fg_0)$, and denote by $\rm{U}_q(\fl_{\le
0})$ the Hopf sub-superalgebra of $\Ul$ generated by all the
generators of $\Ul$ but $E_{m, m+1}$. Note that if
$\Ul\subset\U_q(\fg_0)$, then $\rm{U}_q(\fl_{\le 0})=\Ul$. Let
$L_\mu^{(\fl_{\le 0})}$ be the irreducible $\rm{U}_q(\fl_{\le
0})$-module with integral $\U_q(\fl_0)$ highest weight $\mu$. Note
that the generator $E_{m+1, n}$ of $\rm{U}_q(\fl_{\le 0})$
necessarily acts on $L_\mu^{(\fl_{\le 0})}$ by zero. Furthermore,
$L_\mu^{(\fl_{\le 0})}$ restricts to an irreducible
$\U_q(\fl_0)$-module.
\begin{definition}\bea K_\mu^{(\fl)}:= \Hom_{\U_q(\fl_{\le
0})}\left(\Ul, \, L_\mu^{(\fl_{\le 0})} \right). \label{dKac} \eea
\end{definition}
This will be referred to as a dual Kac module over $\Ul$. The
action of any $y\in\Ul$ on $\zeta\in K_\mu^{(\fl)}$ is given by
$\langle y \zeta, \, x\rangle = (-1)^{[y]([x]+[\zeta])} \langle
\zeta, \, x y\rangle$, $\forall x\in\Ul$.

Let $Ker\epsilon_{\le 0}$ be the subspace of $\U_q(\fl_{\le 0})$
annihilated by the co-unit $\epsilon$. It generates a two-sided
ideal $J(\fl)$ of $\Ul$. By using the PBW theorem \ref{PBW}
(generalized in the obvious way to $\Ul$) we can easily show that
\bea K_\mu^{(\fl)} &=&\left(\Ul/J(\fl)\right)^* \otimes_{\Cq}
L_\mu^{(\fl_{\le 0})}.\label{dKac-mod}\eea It again follows from
the PBW theorem that $\Ul/J(\fl)$ is finite dimensional. The
$\Uh$-module structure of its dual vector space can be described
explicitly. Let $\Phi_1^+(\fl)$ be the set of the odd positive
roots of $\fl$. Let $E$ be the $\Uh$-module with a basis
$\{v_\gamma | \gamma \in \Phi_1^+(\fl)\}$ such that $K_a v_\gamma
= q^{-(\gamma, \epsilon_a)} v_\gamma$, for all $a$. The exterior
algebra of $E$ forms a $\Uh$-module, which we denote by
$\Lambda(\fl_{-1})$. Then $\left(\Ul/J(\fl)\right)^*$ is
isomorphic to $\Lambda(\fl_{-1})$. Therefore, as a $\Uh$-module
(with diagonal action), \bea K_\mu^{(\fl)}&\cong&
\Lambda(\fl_{-1}) \otimes_{\Cq} L_\mu^{(\fl_{\le 0})}. \eea  Note
that when $\Ul$ is purely even, $K_\mu^{(\fl)}=L_\mu^{(\fl_{\le
0})}$.

\section{Quantum homogeneous Super Vector Bundles}\label{Bundles}
In this section we introduce quantum homogeneous superspaces
and quantum homogeneous supervector bundles in the context of
the quantum general linear supergroup. Let us
begin be considering the Hopf superalgebra of functions on the
quantum general linear supergroup.
\subsection{Functions on the quantum general linear supergroup}
General references on Hopf (super)algebra are \cite{MM, Mo}. A
treatment of the classical general linear supergroup similar to what to
be presented here is given in \cite{SZ}.
Let $\Ug^*$ denote the $\Z_2$-graded dual vector space of $\Ug$.
It has a natural associative superalgebraic structure induced by
the co-superalgebraic structure of $\Ug$. Denote the
multiplication of $\Ug^*$ by $m_\circ$, then $\langle
m_{\circ}(f\otimes g),\, x \rangle = \langle f\otimes g, \,
\Delta(x)\rangle,$ for all $x\in\Ug$. (We shall use the notations
$\phi(v)$ and $\langle \phi, \, x\rangle$ interchangeably for the
image of $v\in V$ under $\phi\in\Hom_\Cq(V, W)$.)

There exist two gradation preserving $\Ug$-actions on $\Ug^*$, \be
d \tL: \Ug\otimes \Ug^* \rightarrow \Ug^*, &\quad&
      x\otimes f \mapsto d \tL_x (f), \\
d \tR: \Ug\otimes\Ug^* \rightarrow  \Ug^*, &\quad&
      x\otimes f \mapsto d \tR_x (f), \ee
defined by
\be \langle d \tL_x (f), \, y\rangle& =& (-1)^{[x][f]} \langle f, \,
S(x) y\rangle, \\
\langle d \tR_x (f), \, y\rangle &=& (-1)^{[x]} \langle f, \, y
x\rangle, \quad \forall y\in\Ug. \ee It is easy to see that $d
\tL_{x y} = d \tL_x d \tL_y$, and $d \tR_{x y}=d \tR_x d \tR_y$,
for all $x, y\in\Ug$. Straightforward calculations show that each
of these actions converts $\Ug^*$ into a (graded) left
$\Ug$-module. Furthermore, with respect to the module structure
the product map of $\Ug^*$ is a $\Ug$-module homomorphism and the
unit element of $\Ug^*$ is $\Ug$-invariant. Therefore, each of
these actions converts $\Ug^*$ into a left $\Ug$-module
superalgebra. The fact that the product map in $\Ug^*$ is a module
homomorphism means that the operators $d\tR_x$ and $d\tL_x$ behave
as some sort of generalized super derivations. Indeed, for all
$f,\, g\in\Ug^*$, we have \bea d\tR_x(f g) &=& \sum_{(x)}
(-1)^{[x_{(2)}][f]} d\tR_{x_{(1)}}(f) d\tR_{x_{(2)}}(g), \eea
where we have used the standard Sweedler notation for the
co-multiplication of $x$. However, for $d\tL$, we have \bea
d\tL_x(f g) &=& \sum_{(x)} (-1)^{[x'_{(2)}][f]}d\tL_{x'_{(1)}}(f)
d\tL_{x'_{(2)}}(g),\label{modalg}\eea with respect to the opposite
co-multiplication $\Delta'(x)=\sum_{(x)} x'_{(1)}\otimes x'_{(2)}$
of $x$.  The two actions also super-commute in the sense that
$d\tL_x d\tR_y = (-1)^{[x][y]}d\tR_y d\tL_x$, for all $x,
y\in\Ug$.

Let $ \Uc:=\{f\in\Ug^*\, \mid \, \mbox{kernel of $f$ contains a
co-finite ideal of $\Ug$}\}$ be the finite dual \cite{Mo} of
the quantized universal enveloping algebra $\Ug$ of $\fg$. Here we remark again
that  we only consider $\Z_2$-graded subalgebras and (left, right
or two-sided) ideals of $\Ug$ in this paper. We have the following
lemma, which is an adaption of a standard result (see, e.g., Lemma
9.1.1 in \cite{Mo}) on ordinary associative algebras to $\Ug$. In
fact the result is valid for any associative superalgebra.
\begin{lemma}\label{finite}
For any $f\in\Ug^*$, the following conditions are equivalent:
\begin{enumerate}
\item $f$ vanishes on a left ideal of $\Ug$ of finite
co-dimension;
\item $f$ vanishes on a right ideal of $\Ug$ of finite
co-dimension;
\item $f$ vanishes on an ideal of $\Ug$ of finite
co-dimension, thus belongs to $\Uc$;
\item $d \tL_{\Ug}(f)$ is finite dimensional;
\item $d \tR_{\Ug}(f)$ is finite dimensional;
\item $\left(d \tL_{\Ug}\otimes d \tR_{\Ug}\right)(f)$
      is finite dimensional;
\item $m^*(f)\in\Ug^*\otimes \Ug^*$, where $m^*: \Ug^*
\rightarrow \left(\Ug\otimes\Ug\right)^*$ is defined by $\langle
m^*(f),\, x\otimes y\rangle=\langle f,\, x y \rangle$ for all $x,
y\in\Ug$.
\end{enumerate}
\end{lemma}
\begin{proof}
The proof of Lemma 9.1.1 in \cite{Mo} can be extended verbatim to
superalgebras.
\end{proof}
Therefore, $f$ belongs to $\Uc$ if and only if one of the
equivalent conditions are satisfied.  The Lemma in particular
enables us to impose a Hopf superalgebra structure on $\Uc$, with
multiplication $m_{\circ}$, co-multiplication
$\Delta_{\circ}=\left.m^*\right|_{\Uc}$, unit being $\epsilon$,
and co-unit being the unit $\1_{\Ug}$ of $\Ug$.  The antipode
$S_{\circ}$ of $\Uc$ is defined by \be \langle S_{\circ}(f),\,
x\rangle &=& \langle f,\, S(x)\rangle,
 \quad \forall f\in
\Uc, \ x\in \Ug.  \ee Recall that the antipode $S_\circ$ is
invertible since $S$ is. For convenience, we shall drop the
subscript $\circ$ from the notations for all the structure maps
but the antipode of $\Uc$.

A $\Ug$-representation $\pi$ in $d\times d$-matrices is
a superalgebraic map which is homogeneous of
degree $0$. We write $\pi(x)= \left( \pi_{i j}(x)\right)_{i,
j=1}^{d}$ for any $x\in\Ug$. Define $\pi_{i j}\in\Ug^*,$
 $i, j =1, 2, \ldots , d,$ by
\be \langle\pi_{i j}, \, x\rangle:= \pi_{i j}(x), & &\forall i, j,
\ee and call them the matrix elements of $\pi$. Since the kernel
of any finite dimensional representation is an ideal of $\Ug$ with
finite co-dimension,  all the matrix elements of the
representation belong to $\Uc$. Conversely, $\Uc$ is spanned by
the matrix elements of all the finite dimensional representations
of $\Ug$. To see this, we only need to consider an arbitrary
non-zero element $f\in\Uc$. Let $K$ be a graded co-finite ideal of
$\Ug$ contained in the kernel of $f$. Then $U(\fg)/K$ forms a left
$\Ug$-module under left multiplication,
\be \Ug\otimes U(\fg)/K &\rightarrow&  U(\fg)/K, \\
       y\otimes (x+K) &\mapsto& y x+K.
\ee Let $\{x_i + K\}$ be a basis of $\Ug/K$, and denote by $f_{i
j}$ the matrix elements of the associated representation relative
to this basis. Choose $c_i\in\Cq$ such that $\1_{\Ug}+K=\sum_i c_i
x_i +K$, where $\1_{\Ug}+K$ is not contained in the kernel of $f$
as a set since $f\ne 0$. Then $f=\sum_{i, j} c_i \langle f, \,
x_j\rangle  f_{j i}$.
\begin{definition}
Let $\Ag$ be the $\Z_2$-graded subspace of $\Uc$ spanned by the
matrix elements of the $\Ug$-representations furnished by finite
dimensional objects of $\cC(\fg, \fg)$.
\end{definition}
\begin{lemma}
$\Ag$ is a Hopf sub-superalgebra of $\Uc$.
\end{lemma}
\begin{proof}
The space spanned by the matrix elements of any finite dimensional
representation of $\Ug$ is a sub co-algebra of $\Uc$. Thus $\Ag$
forms a sub co-algebra of $\Uc$. Since $\cC(\fg, \fg)$ is closed
under tensor product, $\Ag$ is closed under multiplication. Also,
if a finite dimensional $\Ug$-module belongs to $\cC(\fg, \fg)$,
then so is also its dual. Thus $\Ag$ is stable under the antipode
of $\Uc$.
\end{proof}
\begin{remark}\label{RUint}
From the discussion on matrix elements of finite dimensional
representations we see that $f\in\Uc$ belongs to $\Ag$ if it
satisfies either $d\tL_{K_a}(f) = q^{(\mu, \epsilon_a)} f$,
$\forall a\in{\bf I}$, or $d\tR_{K_a}(f) = q^{(\mu, \epsilon_a)}
f$, $\forall a\in{\bf I}$, for some $\mu\in\fh_\Z^*$.
\end{remark}

\subsection{Quantum homogeneous super vector bundles}
Recall that in classical geometry, a compact manifold can be
recovered from its algebra of continuous functions by the
Gelfand-Naimark theorem. Also, the Serre-Swan theorem establishes
a one to one correspondence between the spaces of the continuous
sections of vector bundles over a compact manifold and the finite
type projective modules of the algebra of continuous functions on
the manifold. These results are taken as the starting point for
non-commutative geometry \cite{Co}, where `manifolds' are replaced
by non-commutative algebras, and `vector bundles' by finitely
generated projective modules. The quantum homogeneous superspaces
and quantum homogeneous supervector bundles to be studied here
are defined in this spirit.

As is well known, all holomorphic functions on a compact
complex manifold are constants. Therefore, the algebra of holomorphic
functions contains little information about the manifold itself.
This problem persists in classical supergeometry \cite{Ma97} and also
quantum geometry \cite{GZ}. However, as shown in \cite{GZ} in the
context of ordinary quantum groups, we can get around the problem by
working with the quantum analogues of smooth functions in a real setting.
To do this, we need some basic notions about $\ast$-Hopf superalgebras
\cite{Z98, ZZ}.

A $\ast$-superalgebraic structure on an associative superalgebra $A$ over
$\Cq$ is a conjugate linear anti-involution $\theta: A\rightarrow A$: for all
$x, y\in A$, $c,c'\in\Cq$, \bea \theta(c x + d y)={\bar c}\theta(x) +
{\bar c'} \theta(y),\quad& \theta(x y)=\theta(y)\theta(x), \quad&
\theta^2(x)=x. \label{astalg}\eea
Here $\bar c$ and $\bar c'$ are defined in the following way. Let $P$
be a complex polynomial in $q$. Then $\bar P$ is the polynomial obtained
by replacing all the coefficients of $P$ by their complex conjugates.
Now if there is another polynomial $Q$ in $q$ such that $c=P/Q$, then
$\bar c= {\bar P}/{\bar Q}$.
\begin{remark}
Note that the second equation in \eqref{astalg} does not involve any
sign factors as one would normally expect of superalgebras.
\end{remark}

We shall sometimes use the notation $(A, \theta)$ for the $\ast$-superalgebra
$A$ with the $\ast$-structure $\theta$.  Let $(B, \theta_1)$ be
another associative $\ast$-superalgebra. Now $A\otimes B$ has a
natural superalgebra structure, with the multiplication defined for
any $a, a'\in A$ and $b, b'\in B$ by \be (a\otimes b) (a'\otimes b')
&=& (-1)^{[b][a']}a a'\otimes b b'. \ee
Furthermore, the following conjugate linear map \bea
\theta\star\theta_1: a\otimes b \mapsto (1\otimes \theta_1(b))
(\theta(a)\otimes 1) = (-1)^{[a][b]} \theta(a)\otimes\theta_1(b)
\label{star}\eea defines a $\ast$-superalgebraic structure on
$A\otimes B$.

Let us assume that $A$ is a Hopf superalgebra with co-multiplication
$\Delta$, co-unit $\epsilon$ and antipode $S$.  If the
$\ast$-superalgebraic structure $\theta$ satisfies  \be
(\theta\star\theta)\Delta=\Delta\theta, &\quad&
\theta\epsilon=\epsilon\theta,  \ee then $A$ is called a Hopf
$\ast$-superalgebra.  Now $$\sigma:=S \theta$$ satisfies
$\sigma^2=id_A$, as follows from the definition of the antipode.

Let $A^0$ denote the finite dual of $A$, which has a natural Hopf
superalgebraic structure. If $A$ is a Hopf $\ast$-superalgebra
with the $\ast$-structure $\theta$, then $\sigma=S
\theta$ induces a map $\omega: A^0\rightarrow A^0$ defined for any
$f\in A^0$ by \bea \langle \omega(f), x\rangle = \overline{\langle f,
\sigma(x)\rangle }, \quad \forall x\in A.  \label{ast}\eea
As can be easily shown \cite{ZZ}, this map $\omega$ gives rise to a
Hopf $\ast$-superalgebraic structure on $A^0$.

In the case of $\Ug$, the following conjugate anti-involution
defines a Hopf $\ast$-superalgebra structure:
\be \theta:  & E_{a, a+1} &\mapsto E_{a+1, a} K_a K_{a+1}^{-1}, \\
&E_{a+1, a} &\mapsto K_a^{-1} K_{a+1} E_{a, a+1}, \quad\forall a\in{\bf I}', \\
& K_b &\mapsto K_b,\quad \forall b\in{\bf I}.  \ee
The classical counter part of this map determines a compact real
form of the complex general linear supergroup. Let $\sigma=S \theta$,
and define
\bea \URg &:=& \{x\in\Ug | \sigma(x)=x\}.\eea
Clearly $\URg$ forms an associative superalgebra over $\Rq$, even though
it may not have a Hopf superalgebra structure. We shall
refer to it as a real form of $\Ug$, as $\Ug=\Cq\otimes_{\Rq}\URg$.

Now $\Uc$ acquires a Hopf $\ast$-superalgebra structure $\omega$
which is induced by $\sigma$. It is easy to show that the image under
$\omega$ of any matrix element of a finite dimensional
$\Ug$-representation with integral weights must again be a matrix element
of a $\Ug$-representation with the same properties. Thus
\begin{lemma}
$\Ag$ forms a Hopf $\ast$-superalgebra.
\end{lemma}
Therefore $\Ag$ should be considered as some `complexification' of the
superalgebra of functions on a compact form of the quantum general
linear supergroup.

Let us denote by $d L_x$ and $d R_x$ respectively the restrictions
of $d\tL_x$ and $d\tR_x$ to $\Ag$. The following definition will
be important for the remainder of the paper. Let $\Ul$ be a reductive
quantum sub-superalgebra of $\Ug$. Set $\URl:= \Ul\cap\URg$.
We now consider the sub-superalgebra of $\Ag$ invariant under the
left translation of $\URl$.
\begin{definition}  Define \bea \cA(\fg, \fl):=
\left\{ f\in \Ag\, \mid \, dL_x (f) = \epsilon(x) f, \ \forall
x\in\URl \right\},\eea where $\epsilon$ is the co-unit of $\Ug$.
\end{definition}
We shall show presently that $\cA(\fg, \fl)$ forms a superalgebra.
In the philosophy of non-commutative geometry
\cite{Co}, $\cA(\fg, \fl)$ defines some virtue
quantum homogeneous superspace.

It is useful to compare the situation with classical supergeometry \cite{Ma97}.
Let $\Lambda$ be a finite dimensional Grassmann
algebra. Take a parabolic subgroup $P$ of $GL(m|n; \Lambda)$ with
Lie superalgebra $\fp=\fl\oplus\fu$ where $\fu$ is some nilpotent
ideal of $\fp$. Let $U(m|n)$ be a compact real form  of $GL(m|n;
\Lambda)$, and set $K=P\cap U(m|n)$. Then we have the symmetric
superspace $U(m|n)/K$. The tensor product of $\Lambda$ with the
classical analogue of $\cA(\fg, \fl)$ should capture the essential
information of the complexification of the superalgebra  of
functions on $U(m|n)/K$.

Since $\Cq\otimes_{\Rq}\URl= \Ul$, we can show that an element
$f$ of $\Ag$ belongs to $\cA(\fg, \fl)$ if and only if
\be dL_x (f) = \epsilon(x) f, \quad \forall
x\in\Ul. \ee
Also, by Remark \ref{RUint}, an element $g$ of $\Uc$ belongs to
$\Ag$ if $d L_k(f)=\epsilon(k) f$, $\forall k\in\Uh.$ Combining these
observations, we arrive at the following result.
\begin{lemma}\label{inv-Uc}
\be \cA(\fg, \fl)&=&
\left\{ f\in \Uc\, \mid \, dL_x (f) = \epsilon(x) f, \ \forall
x\in\Ul \right\}.\ee
\end{lemma}
Using the left $\Ug$-module algebra structure of $\Ag$, we
immediately show that
\begin{lemma}
The $\cA(\fg, \fl)$ is a sub-superalgebra of $\Ag$.
\end{lemma}
\begin{proof} $\Ul$, being a Hopf sub-superalgebra of $\Ug$,
satisfies $\Delta(\Ul)\subset \Ul\otimes\Ul$. If $f,\ g\in\cA(\fg,
\fl)$, then by \eqref{modalg} we have $dL_x(f g)=\epsilon(x) fg$,
$\forall x\in\Ul$. Therefore, $fg\in\cA(\fg, \fl)$.
\end{proof}

Since $\cA(\fg, \fl)$ is non-commutative, there is a distinction
between left and right $\cA(\fg, \fl)$-modules. However,  the two
sides of the story are `mirror images' of each other, thus we
shall consider $\Z_2$-graded left $\cA(\fg, \fl)$-modules only. A
finitely generated projective module over the superalgebra
$\cA(\fg, \fl)$ will be regarded as the space of sections of a
quantum supervector bundle over the quantum homogeneous
superspace.
\begin{definition}\label{bundle}
Let $\Xi$ be a finite dimensional $\Ul$-module, which naturally
restricts to a $\URl$-module. Define
\be\cS(\Xi)&:=&\left\{\zeta\in \Xi\otimes\Ag  \,\mid\, (\id
\otimes d L_x) \zeta = (S(x)\otimes\id)\zeta, \ \forall x\in\URl
\right\}.\ee
\end{definition} Again it can be easily shown that
\bea
\cS(\Xi)&=&\left\{\zeta\in \Xi\otimes\Ag  \,\mid\, (\id
\otimes d L_x) \zeta = (S(x)\otimes\id)\zeta, \ \forall x\in\Ul
\right\}.
\eea
This fact will be used to prove the following result.
\begin{proposition}\label{S-functor}
\begin{enumerate}
\item$\cS(\Xi)$ forms a left
$\cA(\fg, \fl)$-module under the action
\bea \cA(\fg, \fl)\otimes \cS(\Xi) \rightarrow \cS(\Xi), &\quad&
    f\otimes\zeta \mapsto f\zeta, \eea
defined by $ f\zeta :=\sum (-1)^{[f][w_i]}w_i\otimes  f a_i$ for
$\zeta=\sum w_i\otimes a_i$.
\item Every $\Ul$-module map $\phi: \Xi\rightarrow\Xi'$ induces an
$\cA(\fg, \fl)$-module homomorphism
\bea \cS(\phi)=\phi\otimes\id: \cS(\Xi)\rightarrow\cS(\Xi'). \eea
\end{enumerate}
\end{proposition}
\begin{proof}
For $f\in\cA(\fg, \fl)$ and $\zeta=\sum w_i\otimes
a_i\in\cS(\Xi)$, we have \be (\id \otimes d L_x) f\zeta &=& \sum
(-1)^{([f]+[x])[w_i]}w_i\otimes d L_x (f a_i) \\
&=&\sum (-1)^{([f]+[x])[w_i]+[x][f]}w_i\otimes f \, d L_x (a_i)\\
&=& (S(x)\otimes\id)f\zeta,  \quad \forall x\in\Ul,  \ee where
the second step uses \eqref{modalg} and the defining property of
$\cA(\fg, \fl)$, while the third uses the definition of
$\cS(\Xi)$. This proves the first claim.

The second claim is quite obvious.
\end{proof}

A quantum homogeneous supervector bundle is called trivial if it is
isomorphic to a free left $\cA(\fg, \fl)$-module. Quantum homogeneous
supervector bundles induced by $\Ug$-modules are all trivial.
\begin{proposition}  $\cS(\Xi)$ is freely
generated over $\cA(\fg, \fl)$ if $\Xi$ is the restriction of a
finite dimensional $\Ug$-module.
\end{proposition}
\begin{proof} The proof is adapted from \cite{GZ}.
Note that a finite dimensional left $\Ug$-module
$\Xi$ has a natural right $\Uc$-comodule structure
\be \delta:  \Xi \rightarrow \Xi\otimes \Uc, &
\delta(v)(x)= (-1)^{[v][x]}x v, &\forall v\in\Xi, \ x\in\Ug.\ee
Let $p: \Uc\otimes\Uc\rightarrow \Uc\otimes\Uc$ be
defined by $f\otimes g\mapsto (-1)^{[f][g]}g\otimes f$. Define
a map $\kappa: \Xi\otimes\Uc \rightarrow \Xi\otimes\Uc$ by the composition
of the following  maps
\be \Xi\otimes\Uc\stackrel{\delta\otimes\id}{\longrightarrow}\Xi\otimes\Uc\otimes\Uc
\stackrel{\id\otimes p(S^{-1}\otimes\id)}{\longrightarrow}\Xi\otimes\Uc\otimes\Uc
\stackrel{\id\otimes m_o}{\longrightarrow}\Xi\otimes\Uc, \ee
where $m_o$ is the multiplication of $\Uc$. Explicitly, \bea
\zeta&=&\sum v^{(i)}\otimes f^{(i)} \in \Xi\otimes\Uc,\nonumber\\
 \kappa(\zeta) &=&\sum(-1)^{[f^{(i)}][v_{(2)}^{(i)}]}
v_{(1)}^{(i)}\otimes f^{(i)} S^{-1}(v_{(2)}^{(i)}),\label{kappa}
\eea where we have used Sweedler's notation for $\delta(v^{(i)})$.
The inverse of
$\kappa$ is given by the composition of the following  maps
\be \Xi\otimes\Uc\stackrel{\delta\otimes\id}{\longrightarrow}\Xi\otimes\Uc\otimes\Uc
\stackrel{\id\otimes p}{\longrightarrow}\Xi\otimes\Uc\otimes\Uc
\stackrel{\id\otimes m_o}{\longrightarrow}\Xi\otimes\Uc.  \ee
The restriction of $\kappa$ to $\cS(\Xi)$ is a left
$\cA(\fg, \fl)$-module map as can be easily seen from \eqref{kappa}.
By using \eqref{kappa}, we can also show by a direct calculation that
for any $\zeta\in \cS(\Xi)$,
\be (\id\otimes d L_u)\kappa(\zeta) &=& \epsilon(u) \zeta, \quad \forall
u\in\Ul, \ee
that is $\kappa(\zeta)\in\Xi\otimes \cA(\fg, \fl)$ by Lemma \ref{inv-Uc}. Since $\kappa$ is
invertible, $\kappa\left(\cS(\Xi)\right)=\Xi\otimes \cA(\fg, \fl)$.
\end{proof}

As an immediate consequence of the proposition, we obtain the
following sufficient condition which renders $\cS(\Xi)$ projective
over $\cA(\fg, \fl)$.
\begin{corollary} \label{trivial} The $\cS(\Xi)$ is
projective over $\cA(\fg, \fl)$ if there exists a $\Ul$-module
$\Xi^\bot$ and a finite dimensional $\Ug$-module $V$ such that
$\Xi\oplus\Xi^\bot$ is isomorphic to the restriction of $V$ to a
$\Ul$-module.
\end{corollary}

If $\Ug$ was an ordinary quantized universal enveloping algebra
associated with a finite dimensional semi-simple Lie algebra, it
was shown in \cite{GZ} that $\cS(\Xi)$ was always projective over
$\cA(\fg, \fl)$. Unfortunately this is no longer true for
$\U_q({\mathfrak{gl}}_{m|n})$. However, if $\Ul$ is a purely even
reductive quantum subalgebra of $\Ug$, that is,
$\Ul\subset\rm{U}_q(\fg_0)$, then $\cS(\Xi)$ is a finitely
generated projective $\cA(\fg, \fl)$-module. More generally, we
have the following result. \begin{lemma} If $\lambda\in\fh_\Z^*$
is $\fg$-dominant, then $\cS(K_\lambda^{(\fl)})$ is projective
over $\cA(\fg, \fl)$, where $K_\lambda^{(\fl)}$ is the dual Kac
module over $\Ul$ defined by \eqref{dKac}.\end{lemma}
\begin{proof} To prove this, we let ${\bar
L}_{-\lambda}^{(\fg_{\le 0})}$ be the irreducible $\U_q(\fg_{\le
0})$-module with lowest weight $-\lambda$, which is finite
dimensional. Let ${\bar V}_{-\lambda}= \Ug\otimes_{\U_q(\fg_{\le
0})} {\bar L}_{-\lambda}^{(\fg_{\le 0})}$.  Then ${\bar
V}_{-\lambda}$ is a finite dimensional $\Ug$-module, which
naturally restricts to a $\Ul$-module. Let ${\bar v}_{-\lambda}\in
{\bar V}_{-\lambda}$ be a non-zero vector with weight $-\lambda$, which
generates a $\Ul$-module ${\bar K}=\Ul{\bar v}_{-\lambda}$.
Regard  ${\bar V}_{-\lambda}$ as a $\Uh$-module, we have the decomposition
${\bar V}_{-\lambda}={\bar K}\oplus {\bar K}^\bot$. This in fact is
also a direct sum of $\Ul$-modules as the weights of ${\bar K}^\bot$
differ from those of ${\bar K}$ by roots not belonging to $\fl$.
The dual ${\bar V}^*_{-\lambda}$ of
${\bar V}_{-\lambda}$ has a natural $\Ul$-module structure, and
contains the $\Ul$-submodule ${\bar K}^*=K_\lambda^{(\fl)}$ as a
direct summand.  Therefore Proposition \ref{trivial} applies to
the present situation.\end{proof}

The space $\cS(\Xi)$ has a direct bearing on the representation
theory of $\Ug$.
\begin{lemma}
\begin{enumerate}
\item $\cS(\Xi)$ forms a $\Ug$-module under the action \bea \Ug\otimes
\cS(\Xi) \rightarrow \cS(\Xi), &\quad&
        x \otimes \zeta \mapsto (\id\otimes d R_x)\zeta.
\eea
\item for every $\Ul$-module map
$\phi: \Xi\rightarrow\Xi'$,
the induced map
\bea \cS(\phi)=\phi\otimes\id: \cS(\Xi)\rightarrow\cS(\Xi'), \eea
is $\Ug$-equivariant.
\end{enumerate}
\end{lemma}
\begin{proof}
The first part follows from the fact that the two actions $d L$
and $d R$ of $\Ug$ on $\Ag$ super-commute. To see the second part,
let $\zeta=\sum v_i\otimes f_i$ be in $\cS(\Xi)$. Then for all
$x\in\Ug$,
\be x\circ (\cS(\phi)\zeta) &=& \sum (-1)^{([v_i]+[\phi])[x]}
\phi(v_i)\otimes d R_x (f_i)\\
&=& (-1)^{[x][\phi]} \cS(\phi)(x\circ\zeta).\ee
\end{proof}

Of particular interest to us is the case when $\Xi$ is a finite
dimensional $\Up$-module, where $\Up$ is a parabolic quantum
sub-superalgebra of $\Ug$ with $\Ul$ as its Levi factor. Then
$\cS(\Xi)$ contains the following subspace.
\begin{definition}
$\Gamma(\Xi):=\left\{\zeta\in  \cS(\Xi) \, \mid \, (\id\otimes d
L_x) \zeta =(S(x)\otimes \id)\zeta, \ \forall x\in\Up \right\}.$
\end{definition}
Again by using the super-commutativity of the $\Ug$-actions $d L$
and $d R$ on $\Ag$ we can easily show that
\begin{lemma}
$\Gamma(\Xi)$ is a $\Ug$-submodule of $\cS(\Xi)$. Also a
$\Up$-module homomorphism $\phi: \Xi\rightarrow\Xi'$ induces
a $\Ug$-equivariant map
\bea \Gamma(\phi)=\phi\otimes\id: \Gamma(\Xi)\rightarrow\Gamma(\Xi'). \eea
\end{lemma}

Let $X(\fp, \fl)$ denote the set of $E_{a+1, a}$ or $E_{a, a+1}$
which are contained in $\Up$ but not in $\Ul$. If the $\Up$-module
$\Xi$ has the property that every  element of $X(\fp, \fl)$ acts
by zero, then in this case the definition of $\Gamma(\Xi)$ reduces
to
$$\Gamma(\Xi)= \left\{\zeta\in \cS(\Xi) \, \mid \, (\id\otimes d
L_x ) \zeta =0, \ \forall x\in X(\fp, \fl)\right\}.$$ Thus
$\Gamma(\Xi)$ plays a similar role as the space of holomorphic sections in
classical geometry. We shall refer to it as the space of holomorphic sections
of the homogeneous supervector bundle determined by $\cS(\Xi)$.

We shall promote $\Gamma$ to a covariant functor from the category
$\cC(\fp, \fl)$ of the $\Ul$-finite modules over the parabolic
subalgebra $\Up$ to the category $\cC(\fg, \fg)$ of locally finite
$\Ug$-modules. The resultant functor is shown to be left exact,
and its right derived functors will be regarded as the Dolbeault
cohomology groups of the homogeneous supervector bundle.

\section{Induction Functors}\label{Functors}
We study induction functors and their derived functors in this
section. Results will be applied to develop a representation
theoretical formulation of a quantum analogue of Dolbeault
cohomology for the quantum homogeneous supervector bundles.
References for background material on this section are \cite{EW}
and \cite{KV}. The elementary facts
from homological algebra used here can be found in
any text book, e.g., \cite{We}.
\subsection{Categories of modules}
We start with a discussion on module categories of $\Ug$
and its quantum sub-superalgebras which will be used later.

For any quantum superalgebra $\U$, we shall assume that every
$\U$-module to be considered in this paper is $\Z_2$-graded. Thus
corresponding to each $\U$-module $V$, there exists another module
$\wp V$ which is equal to $V$ as a set, but with $(\wp V)_{\bar
0}= V_{\bar 1}$, and $(\wp V)_{\bar 1}= V_{\bar 0}$. Let $\wp(v)$
denote the element of $\wp V$ corresponding to $v\in V$. The
action of $\U$ on $\wp V$ is defined by $x \wp(v) = (-1)^{[x]}\wp(
x v)$, for all $x\in \U$.

Let $\Ur$ be a quantum sub-superalgebra of $\Ug$. Denote by
$\cC_{inh}(\fr)$ the category of $\Ur$-modules, where the space
$\Hom_{\U}(V, W)$ of morphisms between any two $\U$-modules $V$
and $W$ is a $\Z_2$-graded subspace of $\Hom_{\Cq}(V, W)$
consisting of such elements $\phi$ that for all $x\in \Ur$ and $
v\in V$, $\phi(x v) = (-1)^{[x][\phi]} x \phi(v)$. The parity
change map $\wp$ is odd, and becomes a covariant functor on
$\cC_{inh}(\fr)$ if for any $\phi\in\Hom_{\Ur}(V, W)$ we define
${\wp}(\phi)\in\Hom_{\Ur}(\wp V, \wp W)$ to be the same as $\phi$
on sets.  Note that if $\phi\in \Hom_{\Ur}(V, W)$ is an
inhomogeneous morphism between objects $V$ and $W$ in
$\cC_{inh}(\fr)$, the kernel and image of $\phi$ are not
necessarily $\Z_2$-graded in general, thus
$\cC_{inh}(\fr)$ is not an Abelian category.

Assume $\Ur$ contains a reductive sub-superalgebra $\Uk$ of $\Ug$.
Every $\Ur$-module $V$ naturally restricts to a $\Uk$-module.
\begin{definition}\label{fk-finite} The $\Uk$-finite subspace $V[\Uk]$
of $V$ is defined to be the $\Cq$-span of the integral weight
vectors $v\in V$ satisfying  $\dim(\Uk v )<\infty$.
\end{definition} Here $\Uk v:=\{x v \, \mid \, x\in\Uk\}.$
Elements of $V[\Uk]$ will be called $\Uk$-finite. Also, a
$\Ur$-module $V$ is called $\Uk$-finite if $V=V[\Uk]$.
\begin{remark}\label{Fimage}
If $V$ is a $\Z_2$-graded $\Uk$-finite $\Ur$-module and
$\phi\in\Hom_{\Ur}(V, W)$ a {\em homogeneous} morphism, then
$\phi(V)$ is a $\Z_2$-graded $\Uk$-finite $\Ur$-submodule of $W$.
\end{remark}

Let $\Uq$ be either a parabolic or reductive quantum
sub-superalgebra of $\Ug$. If $\Uq$ contains the reductive quantum
sub-superalgebra $\Uk$, we shall talk about the pair $(\Uq, \Uk)$
of quantum sub-superalgebras. Two pairs of sub-superalgebras
$(\Uq, \Uk)$ and $(\Up, \Ul)$ are said to be compatible if we have
the Hopf superalgebra inclusions $\Uq\supseteq \Up$ and
$\Uk\supseteq\Ul$, and in this case, we write $(\Uq,
\Uk)\supseteq(\Up, \Ul)$. Given a pair $(\Uq, \Uk)$, we shall
denote by $\cC_{inh}(\fq, \fk)$ the full subcategory of
$\cC_{inh}(\fq)$ with the $\Uk$-finite $\Uq$-modules as its
objects. Clearly, $\cC_{inh}(\fq, \fk)$ is closed under passage to
graded sub-modules, graded quotients and finite direct sums. It is
also closed under finite tensor products.
\begin{definition}
Let $\cC(\fq)$ be the subcategory of $\cC_{int}(\fq)$ consisting
of the same objects and the {\em even} morphisms of
$\cC_{inh}(\fq)$. Let $\cC(\fq, \fk)$ be the full subcategory of
$\cC(\fq)$ with the $\Uk$-finite objects.
\end{definition}
Then $\cC(\fq)$ is obviously an Abelian category, and it follows
from Remark \ref{Fimage} that $\cC(\fq, \fk)$ is also Abelian.

\subsection{Induction functors}
Let $(\Up, \Ul)$ be a pair of quantum sub-superalgebras of $\Ug$,
where $\Up$ is either a parabolic or reductive quantum
sub-superalgebra of $\Ug$, and $\Ul$ is a reductive quantum
sub-superalgebra of $\Ug$ contained in $\Up$. We recall that all
the objects of the categories $\cC(\fp)$ and $\cC(\fp, \fl)$ are
$\Z_2$-graded, and all the morphisms of the categories are even.
\begin{definition}\label{Z}
Define a covariant functor $Z_{\fp}^{\fp, \fl}:
\cC(\fp)\rightarrow\cC(\fp, \fl)$ in the following way: for any
object $V$,  $Z_{\fp}^{\fp, \fl}(V)$ is the (not necessarily
direct ) sum of the $\Ul$-finite $\Z_2$-graded $\Up$-submodules of
$V$, and for any morphism $\phi\in\Hom_{\cC(\fp)}(V, W)$,
$Z_{\fp}^{\fp, \fl}(\phi) =\left.\phi\right|_{Z_{\fp}^{\fp,
\fl}(V)}$.
\end{definition}
$Z_{\fp}^{\fp, \fl}$ is well defined because of Remark \ref{Fimage}. When
$\Up=\Ul)=\Ug$, we have an analogue of the Zuckerman functor.
\begin{lemma}
The functor  $Z_{\fp}^{\fp, \fl}: \cC(\fp)\rightarrow\cC(\fp,
\fl)$ is left exact.
\end{lemma}
\begin{proof}
Even though the proof is straightforward, we nevertheless give the details
here because of the importance of the functor $Z_{\fp}^{\fp, \fl}$.
Let us temporarily use $Z$ to denote $Z_{\fp}^{\fp, \fl}.$
Given any exact sequence $$0\longrightarrow U
\stackrel{i}{\longrightarrow} V \stackrel{j}{\longrightarrow} W
$$ in $\cC(\fp), $ we want to show that the
following sequence in $\cC(\fp, \fl)$ is also exact: \be
0\longrightarrow Z(U) \stackrel{Z(i)}{\longrightarrow} Z(V)
\stackrel{Z(j)}{\longrightarrow} Z(W). \ee Assume $U'$ is a
$\Ul$-finite $\Up$-submodule of $U$. Then $i(U')$ is a
$\Ul$-finite $\Up$-submodule of $V$. Thus the injectivity of $i$
implies the injectivity of $Z(i)$.

Let $V'$ be a $\Up$-submodule of $Z(V)$. If an element $v\in V'$
is in $Ker Z(j)$, then there exists a unique $u\in U$ such that
$v=i(u)$. Now $\Up v=i\left(\Up u\right)$ is a $\Ul$-finite
$\Up$-submodule of $V$.  The injectivity of $i$ forces $\Up u$ to
be a $\Up$-submodule of $Z(U)$. In particular, $u\in Z(U)$. Thus
$Im Z(i)\supseteq Ker Z(j)$. But it is obvious that $Im
Z(i)\subseteq Ker Z(j)$. Hence the sequence is also exact at
$Z(V)$.
\end{proof}

Let $\Uq$ either be a parabolic or reductive quantum
sub-superalgebra of $\Ug$, and  let $(\Up, \Ul)$ be as given above
with $\Uq\supseteq\Up$. We define the covariant functor $I_{\fp,
\fl}^{\fq}: \cC(\fp, \fl)\rightarrow \cC(\fq)$ by \bea \I_{\fp,
\fl}^{\fq}(V) := \Hom_{\Up}(\Uq, \, V), && \I_{\fp,
\fl}^{\fq}(\phi) := \Hom_{\Up}(\Uq, \, \phi)\eea for any object
$V$ and morphism $\phi\in\Hom_{\cC(\fp,\fl)}(V, \, W)$. Here
$\I_{\fp, \fl}^{\fq}(\phi)$ is defined for any $\zeta\in \I_{\fp,
\fl}^{\fq}(V)$ by \be  \langle \I_{\fp, \fl}^{\fq}(\phi)(\zeta),
\, x\rangle &=& \phi(\langle \zeta, \, x\rangle),  \quad \forall
x\in \Uq. \ee  Note that $\langle \zeta, \, x \rangle\in V$. The
$\Uq$ action on $\I_{\fp, \fl}^{\fq}(V)$ \bea \Uq\otimes \I_{\fp,
\fl}^{\fq}(V) \rightarrow \I_{\fp, \fl}^{\fq}(V),&\quad& y\otimes
\zeta \mapsto y\circ\zeta, \label{action} \eea is defined by
$\langle y \circ \zeta, \, x\rangle = (-1)^{[y]([x]+[\zeta])}
\langle \zeta, \, x y\rangle$, for all $x\in\Uq$.  The functor
$I_{\fp, \fl}^{\fq}$ is the composition of the exact functor of
tensoring with $\Uq^*$ and the left exact functor of taking $\Up$
invariant submodules, thus is also left exact.

\begin{definition}\label{ZI}
Given compatible pairs $(\Uq, \Uk)\supseteq(\Up, \Ul)$, we
introduce the covariant functor \be\I_{\fp, \fl}^{\fq,
\fk}:=Z_{\fq}^{\fq, \fk}\circ \I_{\fp, \fl}^{\fq}: \, \cC(\fp,\fl)
\rightarrow \cC(\fq, \fk),\ee and call it the induction functor
from $\cC(\fp,\fl)$ to $\cC(\fq, \fk)$.
\end{definition}
\begin{lemma}
The induction functor $\I_{\fp, \fl}^{\fq, \fk}$ is left exact.
\end{lemma}
\begin{proof}
Since both $\I_{\fp, \fl}^{\fq}$ and $Z_{\fq}^{\fq, \fk}$ are left
exact, their composition must also be left exact.
\end{proof}

Let us examine some further properties of the Zuckerman functor
and the induction functor.

Let $\Uq$ be either a parabolic or reductive quantum
sub-superalgebra of $\Ug$, and let $\Uk$ be a reductive quantum
sub-superalgebra of $\Ug$. Assume $\Uq\supset\Uk$. Then for any
object $W$ of $\cC(\fq)$, \bea Z_{\fq}^{\fq, \fk}(W) &=& W[\Uk].
\label{ZW} \eea To prove this, we define the adjoint action of
$\Uq$ on itself \be \ad_y(x) &=& \sum_{(y)} (-1)^{[y_{(2)}][x]}
y_{(1)} x S(y_{(2)}), \quad x, y\in \Uq, \ee where Sweedler's
notation $\Delta(y)=\sum_{(y)} y_{(1)}\otimes y_{(2)}$ is used for
the co-multiplication of $y$. By using the Poincar\'{e}-Birkhoff-Witt
(PBW) Theorem \ref{PBW}, we
can choose a set of $y_i$ each of which is a product of $E_{a b}$
associated with the roots of $\fq$ not contained in $\fk$, such
that every element $x\in\Uq$ can be expressed as a finite sum
$x=\sum y_i u_i$ with $u_i\in\Uk$. By considering the PBW theorem
again, we see that there exists a finite set $Y_x$ of the $y_i$
such that every element of the space $\ad_{\Uk}(x):=\{\ad_u(x) \,
\mid \, u \in\Uk\}$ can be expressed in the form $\sum y_i u_i'$
with $y_i\in Y_x$ and $u_i'\in\Uk$. If $w\in W[\Uk]$,
then \be u (x w)&=& \sum_{(u)} (-1)^{[u_{(2)}][x]}
\left(\ad_{u_{(1)}}(x)\right) \left( u_{(2)} w \right), \quad
x\in\Uq, u\in\Uk. \ee
This implies $u(x w) \in \sum_{y\in Y_x} y \left(\Uk(x
w)\right)$, for all $u\in\Uk$. Therefore,
\be \dim\left(\Uk(x
w)\right)\le |Y_x| \dim\left(\Uk w \right) <\infty, \ee where
$Y_x$ is the cardinality of $Y_x$. Also, if $x$ carries a fixed
weight and $w$ is a weight vector of $W[\Uk]$, then $x w$ is a
weight vector of $W[\Uk]$ with integral weight. Hence $W[\Uk]$ is
indeed a $\Uq$-submodule of $W$.

\begin{lemma}\label{composite} Given compatible pairs
$(\Ur, \Uj))\supseteq(\Uq, \Uk)\supseteq(\Up, \Ul)$ of quantum
sub-superalgebras of $\Ug$, we have  $ \I^{\fr, \fj}_{\fq, \fk}
\circ \I^{\fq, \fk}_{\fp, \fl} =\I^{\fr, \fj}_{\fp, \fl}$ as
covariant functors from $\cC(\fp, \fl)$ to $\cC(\fr, \fj)$.
\end{lemma}
\begin{proof}It is clearly true that for any morphism
$\phi\in\Hom_{\cC(\fp, \fl)}(V, V')$, we have $\I^{\fr, \fj}_{\fq,
\fk} \circ \I^{\fq, \fk}_{\fp, \fl}(\phi)$  $=$  $\I^{\fr,
\fj}_{\fp, \fl}(\phi)$. To prove that the Lemma holds on objects,
we need the following technical result which will be proved below:
if $(\Ur, \Uj)\supseteq(\Uq, \Uk)$, then for any object $W$ of
$\cC(\fq)$, \bea Z_{\fr}^{\fr, \fj}\circ\Hom_{\Uq}\left(\Ur, \,
W\right) &=& Z_{\fr}^{\fr, \fj}\circ\Hom_{\Uq}\left(\Ur, \,
Z^{\fq, \fk}_{\fq}(W)\right). \label{surpass}\eea With
(\ref{surpass}) granted, we have  for any object $V$ of $\cC(\fp,
\fl)$ that \be \I^{\fr, \fj}_{\fq, \fk} \circ \I^{\fq, \fk}_{\fp,
\fl}(V)& =& Z^{\fr, \fj}_{\fr}\circ
\Hom_{\Uq}\left(\Ur,\, \I^{\fq, \fk}_{\fp, \fl}(V)\right)\\
&=&Z^{\fr, \fj}_{\fr}\circ\Hom_{\Uq}\left(\Ur,\, \Hom_{\Up}(\Uq,
V) \right).\ee The far right hand side of the equation can be
simplified by using the following relation \be
\Hom_{\Uq}\left(\Ur, \Hom_{\Up}(\Uq, \, V) \right) &=&
\Hom_{\Up}\left(\Uq\otimes_{\Uq} \Ur,\,  V\right)\\
&=&\Hom_{\Up}\left(\Ur,\,   V\right), \ee and we arrive at \be
\I^{\fr, \fj}_{\fq, \fk} \circ \I^{\fq, \fk}_{\fp, \fl}(V)=Z^{\fr,
\fj}_{\fr} \circ\Hom_{\Up}\left(\Ur,\, V\right) =\I^{\fr,
\fj}_{\fp, \fl}(V). \ee

Now we turn to the proof of equation (\ref{surpass}), which is
equivalent to the statement that for any $\zeta\in Z_{\fr}^{\fr,
\fj}\circ\Hom_{\Uq}\left(\Ur, \, W\right)$, \bea \langle \zeta,\
z\rangle \in Z_{\fq}^{\fq, \fk}(W), &\quad& \forall z\in \Ur.
\label{eqsurpass}\eea By the PBW theorem for $\Ur$, there exists a
set of $x_i$, each of which is a product of elements associated
with roots of $\fr$ not contained in $\fq$,  such that every
$z\in\Ur$ can be expressed as a finite sum $\sum y_i x_i$ with
$y_i\in\Uq$. Let $\nu_i$ be the weight of $x_i$, which is a sum of
roots of $\fr$ thus is integral. We have \be \langle u\circ \zeta,
\, x_i\rangle &=& \sum_{(u)}
(-1)^{[\zeta][u_{(2)}]}\pi_{W}(u_{(1)}) \langle \zeta, \,
\ad_{S(u_{(2)})}(x_i)\rangle,    \quad u\in\Uk,\ee where
$u\circ\zeta$ is as defined by (\ref{action}), and $\pi_W$ refers
to the $\Uq$ action on $W$. From this equation we can deduce that
\be \pi_W(u) \langle \zeta, \, x_i\rangle &=& \sum_{(u)}
(-1)^{[u_{(2)}][\zeta]} \langle u_{(1)}\circ\zeta, \,
\ad_{u_{(2)}}(x_i)\rangle, \quad u\in\Uk.  \ee Arguing as in the proof of
\eqref{ZW},  we conclude that for every $x_i$,
there exists a finite set $X_i$ of the $x_j$ such that every
element of $\ad_{\Uk}(x_i)$ can be expressed as $\sum x_j u_j$,
where $x_j\in X_i$ and $u_j\in\Uk$. Now $\langle \zeta, \sum x_j
u_j\rangle =\sum \langle u_j\circ\zeta, \sum x_j\rangle$, and
$\zeta\in Z_{\fr}^{\fr, \fj}\circ\Hom_{\Uq}\left(\Ur, \, W\right)$
implies that the span of $u\circ\zeta$ for all $u\in\Uk$ is finite
dimensional. Therefore, $$ \dim\left(\pi_W(\Uk)\langle\zeta, \,
x_i\rangle\right)<\infty, \ \mbox{that is},\  \langle\zeta, \,
x_i\rangle\in Z_{\fq}^{\fq, \fk}(W), \forall i.$$  Since
$Z_{\fq}^{\fq, \fk}(W)$ is  $\Uq$-stable, we have $\sum_i
\pi_W(\Uq)\langle\zeta, \, x_i\rangle\subset Z_{\fq}^{\fq,
\fk}(W)$. Now every element of $\Ur$ can be expressed as $\sum y_i
x_i$ with $y_i\in \Uq$. By the definition of $Z_{\fr}^{\fr,
\fj}\circ\Hom_{\Uq}\left(\Ur, \, W\right)$, we have $\sum \langle
\zeta, \, y_i x_i\rangle$  $=$ $\sum(-1)^{[y_i][\zeta]}
 \pi_W(y_i) \langle
\zeta, \, x_i\rangle$,  where the right hand side has just been
shown to belong to $Z_{\fq}^{\fq, \fk}(W)$. This proves equation
(\ref{eqsurpass}), thus completes the proof of the Lemma.
\end{proof}

Denote by $\cF^{\fp, \fl}_{\fq, \fk}: \cC(\fq, \fk) \rightarrow
\cC(\fp, \fl)$ the forgetful functor. We shall refer to the next
theorem as Frobenius reciprocity, which in particular implies that
the induction functor $\I_{\fp, \fl}^{\fq, \fk}$ is right adjoint
to the forgetful functor $\cF^{\fp, \fl}_{\fq, \fk}$. The theorem
plays a crucial role in the study of derived functors of the
induction functors.
\begin{theorem}\label{Frobenius}
There exists the natural even isomorphism \bea \Hom_{\Uq}\left(W,
\, \I_{\fp, \fl}^{\fq, \fk}(V)\right)
&\stackrel{\sim}\longrightarrow& \Hom_{\Up}\left(\cF^{\fp,
\fl}_{\fq, \fk}(W), \, V\right),\nonumber\\ \phi &\mapsto&
\phi(\1_{\Uq}), \label{Frob} \eea of $\Z_2$-graded vector spaces
for all $W$ in $\cC(\fq, \fk)$ and $V$ in $\cC(\fp, \fl)$.
\end{theorem}
\begin{proof} Since $W$ is an object of $\cC(\fq, \fk)$,
the image of any vector of $W$ under an arbitrary
$\phi\in\Hom_{\Uq}\left(W, \,  \Hom_{\Up}(\Uq, \, V)\right)$
belongs to $\I_{\fp, \fl}^{\fq, \fk}(V)$. From this we can easily
deduce that $$\Hom_{\Uq}(W, \,  \I_{\fp, \fl}^{\fq, \fk}(V))
\cong\Hom_{\Uq}\left(W, \,  \Hom_{\Up}(\Uq, \, V)\right).$$ The
right hand side can be further rewritten as
$\Hom_{\Up}\left(\Uq\otimes_{\Uq} W, \,  V\right)$. Now
$\Hom_{\Up}\left(\Uq\otimes_{\Uq}W,  \,  V\right)
\cong\Hom_{\Up}\left(\cF^{\fp, \fl}_{\fq, \fk}(W), \, V\right). $
Thus \bea \Hom_{\Uq}\left(W, \, \I_{\fp, \fl}^{\fq,
\fk}(V)\right)&\cong& \Hom_{\Up}\left(\cF^{\fp, \fl}_{\fq,
\fk}(W), \, V\right). \label{isomorphism}\eea This establishes the
claimed isomorphism between the vector spaces. Let us now show
that $\phi\mapsto \phi(\1_{\Uq})$, $\phi\in\Hom_{\Uq}(W, \,
\I_{\fp, \fl}^{\fq, \fk}(V))$, is indeed the required map. The
$\phi(\1_{\Uq})$ (We shall write $1$ for the identity $\1_{\Uq}$
of $\Uq$.) is the evaluation of $\phi$ at the identity of $\Uq$.
Denote by $\circ$ the action of $\Uq$ on $\I^{\fq, \fk}_{\fp,
\fl}(V)$. Then for any $x\in\Uq$ and $w\in W$, we have \be
\phi(1)(x w)&=(-1)^{[x][\phi]}(x\circ\phi)(1)(w) &=\phi(x)(w), \ee
where the symbol $\circ$ refers to the $\Uq$-action on $\I^{\fq,
\fk}_{\fp, \fl}(V)$. If $\phi$ belongs to the  kernel of the map
(\ref{Frob}), then $\phi(x)=0$, $\forall x\in \Uq$. This forces
$\phi=0$. Thus the map (\ref{Frob}) is injective, and because of
the isomorphism (\ref{isomorphism}), it must be a bijection.

We still need to show that $\phi(1)\in\Hom_{\Up} (\cF^{\fp,
\fl}_{\fq, \fk}(W), \, V)$. But this is clear, as the defining
property of $\I_{\fp, \fl}^{\fq, \fk}(V)$ implies \be
\phi(x)(w)&=& (-1)^{[x][\phi]} x \left(\phi(1)(w)\right), \quad
\forall x\in\Up. \ee This completes the proof.
\end{proof}
The following result is an immediate consequence of
Theorem \ref{Frobenius}.
\begin{corollary} \label{preserve}
The induction functor $\I_{\fp, \fl}^{\fq, \fk}$ takes injectives
to injectives.
\end{corollary}
\begin{proof}
If $V$ is an injective object in $\cC(\fp, \fl)$, then
$\Hom_{\Up}\left( \ \cdot\ , V\right)$ is an exact functor from
$\cC(\fp,  \fl)$ to the category of $\Z_2$-graded vector spaces.
Thus by the Frobenius reciprocity of Theorem \ref{Frobenius},
$\Hom_{\Uq}\left(\ \cdot\ , \I_{\fp, \fl}^{\fq, \fk}(V)\right)$ is
exact on $\cC(\fq, \fk)$.  This in turn implies that $\I_{\fp,
\fl}^{\fq, \fk}(V)$ is injective in $\cC(\fq, \fk)$.
\end{proof}

Consider the pair $(\Up, \Ul)$ of quantum sub-superalgebras of
$\Ug$, where $\Ul$ is assumed to be reductive as always.
\begin{corollary} \label{enough}
The category $\cC(\fp, \fl)$ has enough injectives.
\end{corollary}
\begin{proof}
Let $\rm{U}_q(\fl_0)=\Ul\cap\rm{U}_q(\fg_0)$. Every
$\rm{U}_q(\fl_0)$-finite $\Up$-module is also $\Ul$-finite and
vice versa, hence the categories $\cC(\fp, \fl_0)$ and $\cC(\fp,
\fl)$ are identical. Since $\rm{U}_q(\fl_0)$ is the tensor product
of some non-super $\rm{U}_q({\mathfrak{gl}}_k)$'s, every object of
$\cC(\fl_0, \fl_0)$ is semi-simple, thus is injective. Let $V$ be
an $\Ul$-finite $\Up$-module, which can be restricted to an object
of $\cC(\fl_0, \fl_0)$. Now $\I_{\fl_0, \fl_0}^{\fp, \fl}(V)$ is
injective as follows from the above Corollary. By Theorem
\ref{Frobenius} we have the isomorphism \be F: \,
\Hom_{\Up}\left(V, \, \I_{\fl_0, \fl_0}^{\fp, \fl}(V)\right)
&\stackrel{\sim}{\longrightarrow}& \Hom_{\rm{U}_q(\fl_0)}\left(V,
\, V\right). \ee Consider the pre-image of the identity map
$\id_V\in\Hom_{\rm{U}_q(\fl_0)}(V, \, V)$ under $F$, \bea
\iota:=F^{-1}(\id_V): \,  V\rightarrow\I_{\fl_0, \fl_0}^{\fp,
\fl}(V), &\quad& v\mapsto \iota_v, \label{injection} \eea which is
an injective $\Up$-map. It satisfies $\iota_v(x) = (-1)^{[x][v]} x
v$, \ $ \forall x\in\Up$.
\end{proof}
\begin{remark}
The $\I^{\fq, \fk}_{\fp, \fl}$ can be extended to a covariant
functor $\cC_{inh}(\fq, \fk) \rightarrow \cC_{inh}(\fp, \fl)$ in
the obvious way. It also takes injectives to injectives. By using
Frobenius reciprocity, we can also show that the category
$\cC_{inh}(\fp, \fl)$ has enough injectives.
\end{remark}

Now we restrict our attention to the induction functor $\I_{\fp,
\fl}^{\fg, \fg}: \cC(\fp, \fl)\rightarrow \cC(\fg, \fg)$. Since
the Abelian category $\cC(\fp, \fl)$ has enough injectives, and
$\I_{\fp, \fl}^{\fg, \fg}$ is left exact, it makes sense to talk
about its right derived functors \cite{We} $\left(\I_{\fp,
\fl}^{\fg, \fg}\right)^k$, $k\in\Z_+$. We now give a concrete
description of $\left(\I_{\fp, \fl}^{\fg, \fg}\right)^k$. Let $V$
be any object of $\cC(\fp, \fl)$. Then its restriction to a
$\U_q(\fl_0)$-module lies in $\cC(\fl_0, \fl_0)$ and thus is injective.
We construct the following injective resolution of $V$ in
$\cC(\fp, \fl)$, \bea 0{\rightarrow} V
\stackrel{\iota}{\longrightarrow} I^0(V)
\stackrel{\delta^0}{\longrightarrow} I^1(V)
\stackrel{\delta^1}{\longrightarrow}
I^2(V)\stackrel{\delta^2}{\longrightarrow} \cdots
\label{resolution} \eea where the $\Up$-modules and maps are
defined inductively by \bea &I^{k+1}(V):=\I_{\fl_0, \fl_0}^{\fp,
\fl}(I^k(V)/\delta^{k-1}(I^{k-1}(V))), \nonumber\\
&\delta^k:=\iota\circ p: I^k(V) \stackrel{p}{\rightarrow}
I^k(V)/\delta^{k-1}(I^{k-1}(V))
\stackrel{\iota}{\rightarrow}I^{k+1}(V).\eea Here $\iota$ is
similarly defined as in (\ref{injection}), $p$ is the canonical
projection, and \be I^0(V)= \I_{\fl_0, \fl_0}^{\fp, \fl}(V),& &
\delta^{-1}=\iota.\ee The sequence (\ref{resolution}) is clearly a
resolution, with all $I^k(V)$ being injective because of Corollary
\ref{preserve}. We shall call this injective resolution the {\em
standard resolution}. Now apply the left exact covariant functor
$\I_{\fp, \fl}^{\fg, \fg}$ to it and ignore the first term
$\I_{\fp, \fl}^{\fg, \fg}(V)$, we arrive at the following complex
in $\cC(\fg, \fg)$: \bea 0 {\rightarrow} \Omega^0(\fp, \fl; V)
\stackrel{d^0}{\longrightarrow} \Omega^1(\fp, \fl; V)
\stackrel{d^1}{\longrightarrow} \Omega^2(\fp, \fl;
V)\stackrel{d^2}{\longrightarrow} \cdots,  \label{complex}\eea
where \be \Omega^k(\fp, \fl; V):= \I_{\fp, \fl}^{\fg,
\fg}(I^k(V)), \quad d^k :=\left.\Hom_{\Up} (\Ug, \,
\delta^k)\right|_{\I_{\fp, \fl}^{\fg, \fg}(I^k(V))}. \ee

Let us denote by $\Omega(\fp, \fl; V)$ the complex
(\ref{complex}), and by $H^k(\Omega(\fp, \fl; V))$ its cohomology
groups. Then we have \cite{We} $$\left(\I_{\fp, \fl}^{\fg,
\fg}\right)^k(V)= H^k(\Omega(\fp, \fl; V)),\quad k=0, 1, ...,$$
which are independent of the injective resolution
(\ref{resolution}) chosen. Left exactness of the induction functor
implies \be H^0(\Omega(\fp, \fl; V)) &=& \I_{\fp, \fl}^{\fg,
\fg}(V). \ee

\section{Quantum Bott-Borel-Weil Theorem}\label{BBW}
Throughout this section, we shall assume that $(\Up, \Ul)$ is a
pair of quantum sub-superalgebras such that $\Ul$ is a reductive
quantum sub-superalgebra of $\Ug$, and $\Up$ is the parabolic with
$\Ul$ being its Levi factor. For the sake of concreteness, we also
assume that $\Up$ contains the lower triangular Borel subalgebra
$\rm{U}_q(\bar\fb)$ of $\Ug$. We shall also denote
$\rm{U}_q(\fp_0)=\rm{U}_q(\fg_0)\cap\rm{U}_q(\fp_0)$.

The main results of the section are Theorems \ref{main} and
\ref{general}, which might be considered as a form of
Bott-Borel-Weil
theorem for the quantum general linear supergroup.

\subsection{Dolbeault cohomology groups}\label{H}
We first formulate the Dolbeault cohomology groups of the
homogeneous supervector bundles as the derived functor of an
induction functor. First note that the domain of $\cS$ can be
extended to the category $\cC(\fl, \fl)$,  and that of  $\Gamma$
to the category $\cC(\fp, \fl)$. Below we shall consider these more
generally defined $\cS$ and $\Gamma$. We have the following result.
\begin{theorem}\label{SeqI}
$\cS$ coincides with $\I_{\fl, \fl}^{\fg, \fg}$ on $\cC(\fl, \fl)$.
\end{theorem}
\begin{proof} By examining the second part of Proposition \ref{S-functor}, we
easily see that $\cS$ agrees with $\I_{\fl, \fl}^{\fg, \fg}$ on maps.
Let $\Xi$ be an object of $\cC(\fl, \fl)$.
The inclusion $\cS(\Xi)\subseteq \I_{\fl, \fl}^{\fg, \fg}
(\Xi)$ is obvious since $\Xi\otimes \Ag$ is $\Ug$-finite with
respect to the action $\id_\Xi\otimes d R_{\Ug}$.
Any element $\zeta \in \Hom_{\Ul}\left(\Ug,
\, \Xi\right)$ can be expressed as $\zeta=\sum  \xi_i\otimes f_i$,
where $f_i\in\Ug^*$ and $\zeta_i\in\Xi$.
The $\zeta$ belongs to $\I_{\fl, \fl}^{\fg,
\fg}(\Xi)$ only if $\dim(R_{\Ug}(f_i))<\infty$, $\forall i$. This
is equivalent to the condition that all the $f_i$ belong to $\Uc$,
as follows from Lemma \ref{finite}. Since $\Xi$ regarded as a
$\Uh$-module is integral, we may assume that the $\xi_i$ are
weight vectors with integral weights. The defining property of
$\I_{\fl, \fl}^{\fg, \fg}(\Xi)$ requires the $f_i$ be $dL_{\Uh}$
eigenvectors in $\Uc$ with integral weights. Hence by Remark
\ref{RUint}, the $f_i$ must all belong to $\cA(\fg)$.
\end{proof}
In exactly the same manner we can show that
\begin{corollary}\label{GammaeqIp}
$\Gamma(\Xi)=\I_{\fp, \fl}^{\fg, \fg}(\Xi)$ on $\cC(\fp, \fl).$
\end{corollary}
\begin{remark}
In view of the Theorem and this Corollary, we regard the right derived
functors of $\I_{\fp, \fl}^{\fg, \fg}$ as a form of Dolbeault
cohomology of the quantum homogeneous super vector bundles. Thus
we shall use the more suggestive notation $H^{0, k}(G/P, \,
\cS(\Xi))$ to denote $H^k(\Omega(\fp, \fl; \Xi))$.
\end{remark}

\subsection{The computation of cohomology groups}
The rest of the paper is devoted to the computation of
$H^{0, k}(G/P, \, \cS(\Xi))$.
\subsubsection{A special case with $\Ul\subset\rm{U}_q(\fg_0)$}
Denote by \be \cF^{\fg_0, \fg_0}_{\fg_{\le 0}, \fg_0}:
\cC(\fg_{\le 0}, \fg_0) &\rightarrow& \cC(\fg_0, \fg_0), \\
\cF^{\fg_0, \fl}_{\fg_{\le 0}, \fl}: \cC(\fg_{\le 0}, \fl)
&\rightarrow& \cC(\fg_0, \fl),\\ \cF^{\fp_0, \fl}_{\fp, \fl}:
\cC(\fp, \fl) &\rightarrow& \cC(\fp_0, \fl), \ee the forgetful
functors.
\begin{lemma}\label{forget} If  $\Ul\subset\U_q(\fg_0)$, then
we have the following relations:  \bea \cF^{\fg_0,
\fg_0}_{\fg_{\le 0}, \fg_0} \circ \I_{\fg_{\le 0}, \fl}^{\fg_{\le
0}, \fg_0} &=& \I_{\fg_0, \fl}^{\fg_0, \fg_0}\circ \cF^{\fg_0,
\fl}_{\fg_{\le 0}, \fl}, \label{forget1}\\ \cF^{\fg_0,
\fl}_{\fg_{\le 0}, \fl}\circ \I_{\fp, \fl}^{\fg_{\le 0}, \fl} &=&
\I_{\fp_0, \fl}^{\fg_0, \fl}\circ \cF^{\fp_0, \fl}_{\fp, \fl}.
\label{forget2}\eea
\end{lemma}
\begin{proof}
The first relation can be confirmed by directly checking the
functors involved on objects and morphisms. For any object $V$ of
$\cC(\fg_{\le 0}, \fl)$, $\I_{\fg_{\le 0}, \fl}^{\fg_{\le 0},
\fg_0}(V)=V[\rm{U}_q(\fg_0)]$, and thus $\cF^{\fg_0,
\fg_0}_{\fg_{\le 0}, \fg_0} \circ \I_{\fg_{\le 0}, \fl}^{\fg_{\le
0}, \fg_0}(V)$ is $V[\rm{U}_q(\fg_0)]$ regarded as a
$\rm{U}_q(\fg_0)$-module. Also, $\I_{\fg_0, \fl}^{\fg_0,
\fg_0}\circ \cF^{\fg_0, \fl}_{\fg_{\le 0}, \fl}(V)$ $=$
$\cF^{\fg_0, \fl}_{\fg_{\le 0}, \fl}(V)[\rm{U}_q(\fg_0)]$, which
is again $V[\rm{U}_q(\fg_0)]$ regarded as a
$\rm{U}_q(\fg_0)$-module. Now for any $\phi$ $\in$
$\Hom_{\cC(\fg_{\le 0}, \fl)}\left(V, W\right),$ and all $v\in
V[\rm{U}_q(\fg_0)]$, we have $\cF^{\fg_0, \fg_0}_{\fg_{\le 0},
\fg_0} \circ \I_{\fg_{\le 0}, \fl}^{\fg_{\le 0}, \fg_0}(\phi)(v)$
$=$ $\phi(v)$ $=$ $\I_{\fg_0, \fl}^{\fg_0, \fg_0} \circ
\cF^{\fg_0, \fl}_{\fg_{\le 0}, \fl}(v),$ which belongs to
$W[\rm{U}_q(\fg_0)]$.

Now consider the second relation, which obviously holds on
morphisms.  By using the quantum PBW theorem, we can easily show
that $\rm{U}_q(\fg_0)/\rm{U}_q(\fp_0)\cong \rm{U}_q(\fg_{\le 0}
)/\rm{U}_q(\fp)$ under the given conditions on $\Ul$ and $\Up$.
 Therefore, for any object $V$ of $\cC(\fp,
\fl)$, we have the vector space isomorphism \bea
\Hom_{\rm{U}_q(\fp)}\left(\rm{U}_q(\fg_{\le 0}), V \right)&\cong&
\Hom_{\rm{U}_q(\fp_0)}\left(\rm{U}_q(\fg_0), \cF^{\fp_0,
\fl}_{\fp, \fl}(V)\right). \label{HeqH}\eea Denote by $ P:
\Hom_{\Up}\left(\U_q(\fg_{\le 0}), \, V\right) \rightarrow
\Hom_{\U_q(\fp_0)}\left(\U_q(\fg_0), \, V\right)$ the map induced
by the inclusion of $\U_q(\fg_0)$ in $\U_q(\fg_{\le 0})$, \be
\langle P(\zeta), \, x\rangle &=& \langle \zeta, \, x\rangle,\quad
\forall x\in\U_q(\fg_0)\subset\U_q(\fg_{\le 0}).\ee This map is
$\U_q(\fg_0)$-equivariant, as for any $u\in\U_q(\fg_0)$, we have
\be \langle P(u\circ \zeta), \, x\rangle &=
(-1)^{[u]([x]+[\zeta])}\langle \zeta, \, x u \rangle &=\langle
u\circ P(\zeta), \, x\rangle.\ee Now every element in
$\U_q(\fg_{\le 0})$ may be expressed in the form $\sum y_i u_i$
with $y_i\in\U_q(\fp)$ and $u_i\in\U_q(\fg_0)$.   We have $\langle
\zeta, \, \sum y_i u_i\rangle =\sum (-1)^{[y_i][\zeta]}
\pi_V(y_i)\langle P(\zeta), \, u_i\rangle$. Thus $P(\zeta)=0$ if
and only if $\zeta=0$. Therefore, the $\U_q(\fg_0)$-map $P$ is
injective, which must be bijective because of the vector space
isomorphism \eqref{HeqH}.

Since $\Ul\subset\U_q(\fp_0)\subset\U_q(\fg_{0})$, and
$\Ul\subset\Up\subset\U_q(\fg_{\le 0})$, by \eqref{ZW} we have \be
\I_{\fp_0, \fl}^{\fg_0, \fl}\circ\cF^{\fp_0, \fl}_{\fp, \fl}(V)&=&
\Hom_{\rm{U}_q(\fp_0)}\left(\rm{U}_q(\fg_0), \cF^{\fp_0,
\fl}_{\fp, \fl}(V)\right)[\Ul],\\
\I_{\fp, \fl}^{\fg_{\le 0}, \fl}(V)&=&
\Hom_{\rm{U}_q(\fp)}\left(\rm{U}_q(\fg_{\le 0}), V \right)[\Ul].
\ee The restriction of the $\U_q(\fg_{0})$-equivariant map $P$ to
$\I_{\fp, \fl}^{\fg_{\le 0}, \fl}(V)$ now leads to the sought
after $\U_q(\fg_0)$-module isomorphism \bea \cF^{\fg_0,
\fl}_{\fg_{\le 0}, \fl}\circ \I_{\fp, \fl}^{\fg_{\le 0},
\fl}(V)&\cong& \I_{\fp_0, \fl}^{\fg_0, \fl}\circ\cF^{\fp_0,
\fl}_{\fp, \fl}(V). \label{FIeqIF}\eea
\end{proof}
We also have the following easy result.
\begin{lemma} \label{induced}
The functor $\I_{\fg_{\le 0}, \fg_0}^{\fg, \fg}: \cC(\fg_{\le 0},
\fg_0) \rightarrow \cC(\fg, \fg)$ is exact with \be\I_{\fg_{\le
0}, \fg_0}^{\fg, \fg}(V)&=& \Hom_{\U_q(\fg_{\le 0})}\left(\Ug, \,
V\right),\ee for any object $V$ in $\cC(\fg_{\le 0}, \fg_0)$.
\end{lemma}
\begin{proof} We need to show that $\Hom_{\U_q(\fg_{\le 0})}\left(\Ug, \
\cdot \, \right)$ is exact on $\cC(\fg_{\le 0},
\fg_0) $, and  for any $V$ in $\cC(\fg_{\le 0},
\fg_0) $, $\Hom_{\U_q(\fg_{\le 0})}\left(\Ug, \,
V\right)$ is $\U_q(\fg_0)$-finite. It is fairly easy to see that
$\Hom_{\U_q(\fg_{\le 0})}\left(\Ug, \, V\right)$ is spanned by
integral weigh vectors since $V$ is an object of $\cC(\fg_{\le 0},
\fg_0)$. Let ${\mathcal U}^{+1}$ denote the subspace of $\Ug$
spanned by the ordered products of $(E_{i \alpha})^{\theta_{i
\alpha}}$, $i\le m<\alpha$, $\theta_{i \alpha}=0, 1$. Clearly
$\dim {\mathcal U}^{+1} = 2^{m n}$. Then $$\Hom_{\U_q(\fg_{\le
0})}\left(\Ug, \, V\right) \cong \left({\mathcal
U}^{+1}\right)^*\otimes V.$$ This in particular implies that
$\Hom_{\U_q(\fg_{\le 0})}\left(\Ug, \
\cdot \, \right)$ is exact. Given any $x\in\U_q(\fg_0)$ and
$\eta\in{\mathcal U}^{+1}$, there exist $\eta_j\in {\mathcal
U}^{+1}$ and $x_j\in\U_q(\fg_0)$ such that $\eta x=\sum x_j
\eta_j$. Let $\zeta=\sum f_i\otimes v_i$ be in $\left({\mathcal
U}^{+1}\right)^*\otimes V.$ We have \be \langle x\circ\zeta,
\eta\rangle &=& \sum (-1)^{[v_i][\eta]}\langle f_i, \eta x\rangle
v_i\\
&=& \sum (-1)^{[v_i][\eta]}\langle f_i, \eta_j\rangle  \pi_V(x_j)
v_i.\ee Since $V$ is $\U_q(\fg_0)$-finite, we can deduce from this
equation that $\U_q(\fg_0)\circ\zeta$ is finite dimensional for
any $\zeta=\sum f_i\otimes v_i\in\left({\mathcal
U}^{+1}\right)^*\otimes V.$ This completes the proof.
\end{proof}

The following proposition is one of the main results of this
paper. \begin{proposition}\label{key} Let
$\Ul\subseteq\U_q(\fg_0)$ be a reductive quantum subalgebra of
$\Ug$. Let $\Up\supseteq\Ubb$ be the parabolic quantum
sub-superalgebra of $\Ug$ with $\Ul$ as its Levi factor. Let
$L_\lambda^{(\fp)}$ be a finite dimensional irreducible
$\Up$-module with $\Ul$-highest weight $\lambda\in\fh^*_\Z$.
Denote by $L_\lambda^{(\fp_0)}$ the natural restriction of
$L_\lambda^{(\fp)}$  to a $\U_q(\fp_0)$-module. Then \bea H^{0,
k}(G/P, \,  \cS(L_\lambda^{(\fp)}))&=& \Hom_{\U_q(\fg_{\le
0})}(\Ug, \, \left(\I_{\fp_0, \fl}^{\fg_0, \fg_0}\right)^k
(L_\lambda^{(\fp_0)}) ), \eea where $\left(\I_{\fp_0, \fl}^{\fg_0,
\fg_0}\right)^k(L_\lambda^{(\fp_0)})$ is regarded as a
$\U_q(\fg_{\le 0})$-module with $E_{m+1, m}$ acting  by zero.
\end{proposition}
\begin{proof}
By Lemma \ref{composite}, $ \I_{\fp, \fl}^{\fg, \fg} =
\I_{\fg_{\le 0}, \fg_0}^{\fg, \fg} \circ \I_{\fg_{\le 0},
\fl}^{\fg_{\le 0}, \fg_0} \circ \I_{\fp, \fl}^{\fg_{\le 0}, \fl}.
$ By Lemma \ref{forget}, \be  \cF^{\fg_0, \fg_0}_{\fg_{\le 0},
\fg_0} \circ \I_{\fg_{\le 0}, \fl}^{\fg_{\le 0}, \fg_0} \circ
\I_{\fp, \fl}^{\fg_{\le 0}, \fl} =\I_{\fg_0, \fl}^{\fg_0, \fg_0}
\circ \I_{\fp_0, \fl}^{\fg_0, \fl}\circ\cF^{\fp_0, \fl}_{\fp,
\fl}. \ee Using Lemma \ref{composite} again, we obtain \bea
\cF^{\fg_0, \fg_0}_{\fg_{\le 0}, \fg_0} \circ \I^{\fg_{\le 0},
\fg_0}_{\fp, \fl}&=&\I_{\fp_0, \fl}^{\fg_0, \fg_0}
\circ\cF^{\fp_0, \fl}_{\fp, \fl}. \eea Recall the following
elementary facts: Let $\cC\stackrel{G}{\longrightarrow}\cC'$ be a
left exact covariant functor. (a). Suppose
$\cC'\stackrel{F}{\longrightarrow}\cC''$ is an exact covariant
functor. Then $F\circ G$ is left exact, and its right derive
functors are $(F\circ G)^k=F\circ (G)^k$. (b). Suppose
${\tilde\cC}\stackrel{F}{\longrightarrow}\cC'$ is an exact
covariant functor. Then $G\circ F$ is left exact, and its right
derived functors are $(G\circ F)^k=(G)^k\circ F$. Applying these
results to the situation at hand, we arrive at \be  \cF^{\fg_0,
\fg_0}_{\fg_{\le 0}, \fg_0} \circ \left(\I^{\fg_{\le 0},
\fg_0}_{\fp, \fl}\right)^k (L_\lambda^{(\fp)}) &=&
\left(\I_{\fp_0, \fl}^{\fg_0,
\fg_0}\right)^k(L_\lambda^{(\fp_0)}). \ee The derived functor
$\left(\I_{\fp_0, \fl}^{\fg_0, \fg_0}\right)^k$ on the right hand
side can be computed by using the quantum Bott-Borel-Weil theorem
\cite{APW} for $\U_q(\fg_0)=\U_q(\mathfrak{gl}_m)\otimes
\U_q(\mathfrak{gl}_n)$. Now $\left(\I_{\fp_0, \fl}^{\fg_0,
\fg_0}\right)^k(L_\lambda^{(\fp_0)})$ is either zero or a finite
dimensional irreducible $\U_q(\fg_0)$-module. Therefore, its
inverse image under the forgetful functor $\cF^{\fg_0,
\fg_0}_{\fg_{\le 0}, \fg_0}$ must be either zero or $\U_q(\fg_{\le
0})$-irreducible. In both cases, $E_{m+1, m}$ acts by zero.  Thus
by using Lemma \ref{induced}, we have \be \left(\I_{\fp,
\fl}^{\fg, \fg}\right)^k(L_\lambda^{(\fp)}) &=& \I_{\fg_{\le 0},
\fg_0}^{\fg, \fg}\left( \left(\I_{\fp_0, \fl}^{\fg_0,
\fg_0}\right)^k(L_\lambda^{(\fp_0)}) \right),\ee where
$\left(\I_{\fp_0, \fl}^{\fg_0,
\fg_0}\right)^k(L_\lambda^{(\fp_0)})$ is regarded as a
$\U_q(\fg_{\le 0})$-module with $E_{m+1, m}$ acting by zero.
Another easy application of Lemma \ref{induced} completes the
proof.
\end{proof}

By using the proposition we can easily prove the following result.
\begin{theorem} \label{special} \label{main}
Let $\Ul\subseteq\U_q(\fg_0)$ be a reductive quantum subalgebra of
$\Ug$. Let $\Up\supseteq\Ubb$ be the parabolic quantum
sub-superalgebra of $\Ug$ with $\Ul$ as its Levi factor. Let
$L_\lambda^{(\fp)}$ be a finite dimensional irreducible
$\Up$-module with $\Ul$-highest weight $\lambda\in\fh^*_\Z$.
\begin{enumerate}
\item If $\lambda$ is $\fg$-regular, then there exists a unique
element $w$ of the Weyl group of $\fg_0$ rendering
$\mu:=w(\lambda+\rho)-\rho$ dominant with respect to $\fg$. In
this case, \be H^{0, k}(G/P, \, \cS(L_\lambda^{(\fp)}))&=&
\left\{\begin{array}{l l} K_\mu^{(\fg)}, & k=|w|, \\ 0, &k\ne |w|,
\end{array}\right.\ee
where $|w|$ denotes the length of $w$.
\item If $\lambda$ is not $\fg$-regular, then
\be H^{0, k}(G/P, \, \cS(L_\lambda^{(\fp)}))&=& 0,  \quad \forall
k.\ee
\end{enumerate}
\end{theorem}
\begin{proof}
According to the quantum Bott-Borel-Weil theorem for
quantized universal enveloping algebras of ordinary Lie algebras \cite{APW},
the $\left(\I_{\fp_0, \fl}^{\fg_0,
\fg_0}\right)^k (L_\lambda^{(\fp_0)})$ vanishes for all $k$ if
$\lambda$ is not $\fg_0$-regular. If $\lambda$ is $\fg_0$-regular,
then $\left(\I_{\fp_0, \fl}^{\fg_0, \fg_0}\right)^k
(L_\lambda^{(\fp_0)})$ is concentrated at on degree, namely,
$\left(\I_{\fp_0, \fl}^{\fg_0, \fg_0}\right)^k
(L_\lambda^{(\fp_0)})$ is non-vanishing for one $k$ only. We have
\be \left(\I_{\fp_0, \fl}^{\fg_0, \fg_0}\right)^{|w|}
(L_\lambda^{(\fp_0)}) &=& L_\mu^{(\fg_0)}, \ee where
$L_\lambda^{(\fg_0)}$ is the irreducible $\U_q(\fg_0)$-module with
highest weight $\mu$. Using this result in Proposition \ref{key},
we arrive at the Theorem.
\end{proof}
\begin{remark}\label{sum} An easy examination will show that the proof for
Proposition \ref{key} still goes through for
$\U_q\left(\mathfrak{gl}_{m_1|n_1}\oplus
\mathfrak{gl}_{m_2|n_2}\oplus ...\oplus
\mathfrak{gl}_{m_i|n_i}\right)$ for any finite $i$. The same
comment applies to Theorem \ref{main}.
\end{remark}

\subsubsection{The general case} We investigate the general case in
this subsection. Now $\Ul$ is an arbitrary reductive quantum
sub-superalgebra of $\Ug$, and $\Up$ is the parabolic containing
$\Ubb$ and has the Levi factor $\Ul$.  Let
$\U_q(\bar\fb_{\fl})=\Ubb\cap\Ul$ be the Borel subalgebra of
$\Ul$. Denote by \be \cF^{\fl, \fl}_{\fp, \fl}: \cC(\fp, \fl)
\rightarrow \cC(\fl, \fl), &\quad& \cF^{{\bar\fb}_{\fl},
\fh}_{\bar\fb, \fh}: \cC(\bar\fb, \fh) \rightarrow
\cC({\bar\fb}_{\fl}, \fh) \ee the forgetful functors. We have the
following result.
\begin{lemma} \label{forget3}
$\cF^{\fl, \fl}_{\fp, \fl}\circ \I^{\fp, \fl}_{\bar\fb, \fh}
=\I^{\fl, \fl}_{{\bar\fb}_{\fl}, \fh}
\circ\cF^{{\bar\fb}_{\fl}, \fh}_{\bar\fb,
\fh}$.
\end{lemma}
\begin{proof}
The proof is much the same as that for (\ref{forget2}). Because of
the given conditions on $\Up$ and $\Ul$, equation \eqref{ZW} gives
 \be \I^{\fp,
\fl}_{\bar\fb, \fh}(V) &=&\Hom_{\Ubb} \left(\Up, \, V\right)[\Ul],
\ee for any object $V$ of $\cC(\bar\fb, \fh)$. We can easily show
that there exists the even $\Ul$-module isomorphism \be P:  \,
\Hom_{\Ubb}\left(\Up, \, V\right)
&\stackrel{\sim}{\longrightarrow}& \Hom_{\U_q(\bar\fb_{\fl})}
\left(\Ul, \, V\right)\ee defined by $\langle \zeta, \, u x\rangle
= \pi_V(u) \langle P(\zeta), \, x\rangle$, for all $u\in\Ubb$,
$x\in\Ul$. Therefore, \be \I^{\fp, \fl}_{\bar\fb, \fh}(V) &=&
\Hom_{\U_q({\bar\fb}_{\fl})}\left(\Ul, \, V\right)[\Ul]. \ee On
the other hand, \be \I^{\fl, \fl}_{{\bar\fb}_{\fl}, \fh}\circ
\cF^{{\bar\fb}_{\fl}, \fh}_{\bar\fb, \fh}(V)&=&
\Hom_{\U_q({\bar\fb}_{\fl})}\left(\Ul, \, \cF^{{\bar\fb}_{\fl},
\fh}_{\bar\fb, \fh}(V)\right)[\Ul]. \ee Thus the claim of the
Lemma is indeed true for any object of $\cC(\bar\fb, \fh)$. The
claim also clearly holds true for morphisms of $\cC(\bar\fb,
\fh)$.
\end{proof}

\begin{theorem}\label{general}
Let $\lambda\in\fh^*_\Z$ be $\fl$-dominant. Inflate
$K_\lambda^{(\fl)}$ to a $\Up$ module by requiring that all the
generators of $\Up$ not contained in $\Ul$ act by zero, and denote
the resultant $\Up$-module by $K_\lambda^{(\fp)}$.
\begin{enumerate}
\item If $\lambda$ is $\fg$-regular, then there exists a unique $w$
in the Weyl group of $\fg_0$ rendering $\fg$-dominant the
following weight $\mu:=w(\lambda+\rho)-\rho$. In this case,
\be
H^{0, k}(G/P, \,  \cS(K^{(\fp)}_\lambda))&=&
\left\{\begin{array}{l l} K_\mu^{(\fg)}, & k=|w|, \\ 0, &k\ne |w|.
\end{array}\right.\ee
\item If $\lambda$ is not $\fg$-regular, then
$H^{0, k}(G/P, \,  \cS(K^{(\fp)}_\lambda))=0, \quad \forall k.$
\end{enumerate}
\end{theorem}
\begin{proof}
We use Lemma \ref{composite} to write $\I^{\fg, \fg}_{\bar\fb,
\fh} =\I^{\fg, \fg}_{\fp, \fl}\circ \I^{\fp, \fl}_{\bar\fb, \fh}.$
The functor $\I^{\fp, \fl}_{\bar\fb, \fh}$ takes injectives to
injectives. Thus for an irreducible $\Ubb$-module $\Cq_\lambda$
with an arbitrary weight $\lambda\in\fh^*_\Z$, we have a first
quadrant spectral sequence, the Grothendieck spectral sequence
(Sections 5.8  and 10.8 of \cite{We}), \be E_r^{p, q}
&\Longrightarrow & \left(\I^{\fg, \fg}_{\bar\fb,
\fh}\right)^{p+q}(\Cq_\lambda), \ee with $E_2^{p, q}$ term \be
E_2^{p, q} &=& \left(\I^{\fg, \fg}_{\fp, \fl}\right)^{p}
\left(\I^{\fp, \fl}_{\bar\fb, \fh}\right)^{q}(\Cq_\lambda), \ee
where the differential on $E_r^{p, q}$ has bi-degree $(r,\ 1-r)$.
We shall prove below that $\left(\I^{\fp, \fl}_{\bar\fb,
\fh}\right)^{q}(\Cq_\lambda)$ is concentrated at one degree. Let
us take this as granted for the moment.  Then the spectral
sequence collapses at $E_2$, and we obtain \bea  \left(\I^{\fg,
\fg}_{\bar\fb, \fh}\right)^{p+q}(\Cq_\lambda)&=& \left(\I^{\fg,
\fg}_{\fp, \fl}\right)^{p}  \left(\I^{\fp, \fl}_{\bar\fb,
\fh}\right)^{q}(\Cq_\lambda).  \label{spectral}\eea

Now we consider $\left(\I^{\fp, \fl}_{\bar\fb,
\fh}\right)^{q}(\Cq_\lambda)$ for arbitrary
$\lambda\in\fh^*_{\Z}$. By Lemma \ref{forget3}, we have \bea
\cF^{\fl, \fl}_{\fp, \fl}\circ \left(\I^{\fp, \fl}_{\bar\fb,
\fh}\right)^q (\Cq_\lambda) &=& \left(\I^{\fl,
\fl}_{{\bar\fb}_{\fl}, \fh}\right)^q \circ\cF^{{\bar\fb}_{\fl},
\fh}_{\fb, \fh}(\Cq_\lambda). \label{HK} \eea Note that $\Ul$ is
the tensor product of the quantized universal enveloping algebra
of the direct sum of some general linear algebras and possibly
also a general linear superalgebra. By Theorem \ref{main} and
Remark \ref{sum}, the right hand side is zero unless $\lambda$ is
$\fl$-regular. When $\lambda$ is $\fl$-regular, $\left(\I^{\fl,
\fl}_{{\bar\fb}_{\fl}, \fh}\right)^q \circ\cF^{{\bar\fb}_{\fl},
\fh}_{\bar\fb, \fh}(\Cq_\lambda)$ is concentrated at one degree.
Explicitly, there exists a unique $w_\fl$ in the Weyl group of
$\fl$ rendering $w_{\fl}(\lambda+\rho_{\fl})-\rho_{\fl}$ dominant
with respect to $\fl$, and we have \be \left(\I^{\fl,
\fl}_{{\bar\fb}_{\fl}, \fh}\right)^{|w_\fl|}
\circ\cF^{{\bar\fb}_{\fl}, \fh}_{\fb, \fh}(\Cq_\lambda) =
K_{w_\fl(\lambda+\rho_{\fl})-\rho_{\fl}}^{(\fl)}.\ee Here
$\rho_{\fl}$ is half of the signed-sum of the positive roots of
$\fl$ relative to $\fb_{\fl}=\fb\cap\fl$. Needless to say, the
formula remains valid if we replace $\rho_{\fl}$ by $\rho$.

In order to determine $\left(\I^{\fp, \fl}_{\bar\fb, \fh}\right)^q
(\Cq_\lambda)$, we consider all the possible objects $W_\lambda$
of $\cC(\fp, \fl)$ satisfying $\cF^{\fl, \fl}_{\fp,
\fl}(W_\lambda)$ $=$ $K_{w_\fl(\lambda+\rho)-\rho}^{(\fl)}$. Any
two weights of $K_{w_\fl(\lambda+\rho)-\rho}^{(\fl)}$ can only
differ by an integeral combination of the roots of $\fl$. This in
particular requires that all the generators of $\Up$ not contained
in $\Ul$ act on $W_\lambda$ by zero. Therefore, \be \left(\I^{\fp,
\fl}_{\bar\fb, \fh}\right)^{|w_\fl|}(\Cq_\lambda) &=&
K_{w_\fl(\lambda+\rho)-\rho}^{(\fp)}.\ee By using the given
condition that $\lambda$ is $\fl$-dominant, we obtain from
(\ref{spectral}) \be \left(\I^{\fg, \fg}_{\fp,
\fl}\right)^{k}\left(K_{\lambda}^{(\fp)} \right) &=&\left(\I^{\fg,
\fg}_{\bar\fb, \fh}\right)^{k}(\Cq_\lambda). \ee Using the special
case of Theorem \ref{main} with the parabolic being $\Ubb$, we
complete the proof.
\end{proof}

\end{document}